\documentclass[a4paper]{amsart}
\usepackage{graphics}
\usepackage{amsmath}
\usepackage{latexsym}
\usepackage{amssymb}
\usepackage[all]{xy}
\xyoption{matrix}
\xyoption{arrow}

\begin{document}

\renewcommand{\th}{\operatorname{th}\nolimits}
\newcommand{\rej}{\operatorname{rej}\nolimits}
\newcommand{\extto}{\xrightarrow}
\renewcommand{\mod}{\operatorname{mod}\nolimits}
\newcommand{\Sub}{\operatorname{Sub}\nolimits}
\newcommand{\Coh}{\operatorname{Coh}\nolimits}
\newcommand{\ind}{\operatorname{ind}\nolimits}
\newcommand{\Fac}{\operatorname{Fac}\nolimits}
\newcommand{\add}{\operatorname{add}\nolimits}
\newcommand{\red}{\operatorname{red}\nolimits}
\newcommand{\Hom}{\operatorname{Hom}\nolimits}
\newcommand{\Irr}{\operatorname{Irr}\nolimits}
\newcommand{\Rad}{\operatorname{Rad}\nolimits}
\newcommand{\RHom}{\operatorname{\mathbb{R}Hom}\nolimits}
\newcommand{\uHom}{\operatorname{\underline{Hom}}\nolimits}
\newcommand{\End}{\operatorname{End}\nolimits}
\renewcommand{\Im}{\operatorname{Im}\nolimits}
\newcommand{\Ker}{\operatorname{Ker}\nolimits}
\newcommand{\Coker}{\operatorname{Coker}\nolimits}
\newcommand{\Ext}{\operatorname{Ext}\nolimits}
\newcommand{\Tor}{\operatorname{Tor}\nolimits}
\newcommand{\op}{{\operatorname{op}}}
\newcommand{\Ab}{\operatorname{Ab}\nolimits}
\newcommand{\id}{\operatorname{id}\nolimits}
\newcommand{\pd}{\operatorname{pd}\nolimits}
\newcommand{\A}{\operatorname{\mathcal A}\nolimits}
\newcommand{\C}{\operatorname{\mathcal C}\nolimits}
\newcommand{\D}{\operatorname{\mathcal D}\nolimits}
\newcommand{\E}{\operatorname{\mathcal E}\nolimits}
\newcommand{\He}{\operatorname{\mathcal H}\nolimits}
\newcommand{\M}{\operatorname{\mathcal M}\nolimits}
\newcommand{\X}{\operatorname{\mathcal X}\nolimits}
\newcommand{\U}{\operatorname{\mathcal U}\nolimits}
\newcommand{\XX}{\operatorname{\mathbb X}\nolimits}
\renewcommand{\S}{\operatorname{\mathcal S}\nolimits}
\newcommand{\Y}{\operatorname{\mathcal Y}\nolimits}
\newcommand{\F}{\operatorname{\mathcal F}\nolimits}
\newcommand{\Z}{\operatorname{\mathbb Z}\nolimits}
\newcommand{\LL}{\operatorname{\mathbb L}\nolimits}
\newcommand{\Q}{\operatorname{\mathbb Q}\nolimits}
\renewcommand{\P}{\operatorname{\mathcal P}\nolimits}
\newcommand{\T}{\mathcal T}
\newcommand{\G}{\Gamma}
\renewcommand{\L}{\Lambda}
\renewcommand{\r}{\operatorname{\underline{r}}\nolimits}
\newtheorem{lem}{Lemma}[section]
\newtheorem{prop}[lem]{Proposition}
\newtheorem{cor}[lem]{Corollary}
\newtheorem{thm}[lem]{Theorem}
\newtheorem*{thmA}{Theorem A}

\newtheorem*{thmB}{Theorem B}
\newtheorem{rem}[lem]{Remark}
\newtheorem{defin}[lem]{Definition}

\providecommand{\abs}[1]{\lvert#1\rvert}


\title[Cluster mutation via quiver representations]{Cluster mutation via quiver representations}

\author[Buan]{Aslak Bakke Buan}
\address{Institutt for matematiske fag\\
Norges teknisk-naturvitenskapelige universitet\\
N-7491 Trondheim\\
Norway}
\email{aslakb@math.ntnu.no}

\author[Marsh]{Robert J. Marsh}
\address{Department of Mathematics \\
University of Leicester \\
University Road \\
Leicester LE1 7RH \\
England
}
\email{rjm25@mcs.le.ac.uk}

\author[Reiten]{Idun Reiten}
\address{Institutt for matematiske fag\\
Norges teknisk-naturvitenskapelige universitet\\
N-7491 Trondheim\\
Norway}
\email{idunr@math.ntnu.no}


\thanks{Aslak Bakke Buan was supported by a grant from the Norwegian Research Council}

\begin{abstract}
Matrix mutation appears in the definition of cluster algebras of Fomin and Zelevinsky.
We give a representation theoretic interpretation of matrix mutation, using tilting theory
in cluster categories of hereditary algebras. Using this, 
we obtain a representation theoretic interpretation of cluster mutation 
in case of acyclic cluster algebras of finite type.
\end{abstract}

\maketitle

\section*{Introduction}
This paper was motivated by the interplay between the recent development of the theory of 
cluster algebras defined by Fomin and Zelevinsky in \cite{fz1} (see \cite{z} for an introduction)
and the subsequent theory of cluster categories and cluster-tilted algebras \cite{bmrrt, bmr}.
Our main results can be considered to be interpretations within cluster categories of important
concepts in the theory of cluster algebras.

Cluster algebras were introduced in order to explain the
connection between the canonical basis of a quantized enveloping algebra
as defined by Kashiwara and Lusztig and total positivity for algebraic groups.
It was also expected that cluster algebras
should model the classical and quantized coordinate rings of varieties
associated to algebraic groups --- see~\cite{bfz} for an example
of this phenomenon (double Bruhat cells).
Cluster algebras have been used to define generalizations of the Stasheff
polytopes (associahedra) to other Dynkin types~\cite{cfz,fz3}; consequently
there are likely to be interesting links with operad theory.
They have been used to
provide the solution~\cite{fz3} of a conjecture of
Zamolodchikov concerning $Y$-systems, a class of functional relations
important in the theory of the thermodynamic Bethe Ansatz, as well as
solution~\cite{fz4} of various recurrence problems involving
Laurent polynomials, including a conjecture of Gale and Robinson on the
integrality of generalized Somos sequences. Here the
remarkable Laurent properties of the distinguished generators of a cluster
algebra play an important role.
Cluster algebras have also been related to Poisson geometry~\cite{gsv1},
Teichm\"{u}ller spaces~\cite{gsv2}, positive spaces and stacks~\cite{fg1},
dual braid monoids~\cite{bes}, ad-nilpotent ideals of a Borel subalgebra
of a simple Lie algebra~\cite{p} as well as representation
theory, see amongst others~\cite{bmrrt,bmr,cc,ccs1,ccs2,mrz}.

A cluster algebra (without coefficients) is defined via a choice of
free generating set 
$\underline{x} = \{ x_1, \dots, x_n \} $ in the field $\F$ of rational polynomials over $\mathbb{Q}$
and a skew-symmetrizable integer matrix $B$ indexed by
the elements of $\underline{x}$. 
The pair $(\underline{x},B)$, called a seed, determines the cluster algebra as a
subring of $\F$. More specifically, for each
$i= 1, \dots, n$, a new seed $\mu_i(\underline{x},B) = (\underline{x'},B')$ is obtained by replacing
$x_i$ in $\underline{x}$ by 
${x_i}' \in \F$, where ${x_i}'$ is obtained by a so-called
{\em exchange multiplication rule} and
$B'$ is obtained from $B$ by applying so-called
{\em matrix mutation} at row/column $i$.
Mutation in any direction is also defined for the new seed, and by iterating this process
one obtains a countable (sometimes finite) number of seeds. For a seed $(\underline{x},B)$, 
the set $\underline{x}$ is called a {\em cluster}, and the elements in 
$\underline{x}$ are called {\em cluster variables}.
The desired subring of $\F$ is by definition generated by the cluster variables.

It is an interesting problem to try to find a categorical/
module theoretical setting with a nice interpretation of the concepts of clusters and 
cluster variables, and of the matrix mutation and multiplication exchange rule for cluster
variables. For the case of acyclic cluster variables so-called {\em cluster categories} 
were introduced as a candidate for such a model \cite{bmrrt}. Skew-symmetric matrices
are in one-one correspondence with finite quivers with no loops or cycles of length two,
and the corresponding cluster algebra is called {\em acyclic} if there is a seed 
$(\underline{x},B)$ such that $B$ corresponds to a quiver $Q$ without oriented cycles.
There is then, for a field $k$, an associated finite dimensional path algebra $kQ$.
The corresponding cluster category $\C$ is defined in \cite{bmrrt} as a certain quotient 
of the bounded derived category of $kQ$, which is shown to 
be canonically triangulated by \cite{k}. In \cite{bmrrt} (cluster-)tilting theory is developed in
$\C$, with emphasis on connections to cluster algebras. The analogs 
of clusters are (cluster-)tilting objects, and the analogs of cluster variables are
exceptional objects. 
In case $Q$ is a Dynkin quiver, it was shown in \cite{bmrrt} that there is a one-one
correspondence between cluster variables and exceptional objects in $\C$ (in this
case all indecomposables are exceptional) which takes clusters to tilting objects.
This was conjectured to hold more generally.

In this paper we show that also the matrix mutation for cluster algebras has a nice interpretation 
within cluster categories, in terms of the associated cluster-tilted algebras, investigated in
\cite{bmr}. Cluster-tilted algebras are endomorphism algebras of tilting objects in
cluster categories. It follows from our results that the quivers of the cluster-tilted algebras 
arising from a given cluster category are
exactly the quivers corresponding to the exchange matrices of the associated cluster algebra.
This has further applications to cluster algebras (see \cite{br}).
Another main result of this paper is an interpretation within cluster categories
of the exchange multiplication rule
of an (acyclic) cluster algebra. 
So, together with the results from \cite{bmrrt}, all the major ingredients involved in the
construction of acyclic cluster algebras have now been interpreted in the cluster category.

Tilting theory for hereditary algebras have been
a central topic within representation theory since the early eighties. 
This involves the study of tilted algebras, and various generalizations. 
An important motivation for this theory was to compare 
the representation theory of hereditary algebras with the representation theory of
other homologically more complex algebras.
The main result of \cite{bmr} is also in this spirit, showing a close connection between
the representation theory of cluster-tilted algebras and hereditary algebras.
It is the hope of the authors that our ``dictionary'' also can be used
to obtain further developments in representation theory of finite dimensional algebras. 
Also new links between this field and other fields of mathematics can be expected,
having in mind the influence of cluster algebras on other areas.

In \cite{ccs1} an alternative description
of the cluster category is given for type $A$. 
The cluster category was also the motivation for a Hall-algebra type
definition of a cluster algebra of finite type~\cite{cc, ck}.

The paper is organized as follows. In section 1 we give some 
preliminaries, allowing us to state the main result more precisely. Most of the necessary 
background on cluster algebras is however postponed until later (section 6),
since most of the paper does not involve cluster algebras. In section 2 we prove
the following: If $\G$ is cluster-tilted, then so is $\G / \G e \G$ for
an idempotent $e$ in $\G$. This is an 
essential ingredient in the proof of the main result, and also an interesting
fact in itself. In section 3 some consequences of this are given.
In section 4 we prepare for the proof of our main result. This involves
studying cluster-tilted algebras of rank 3, and a crucial result of Kerner \cite{ker}
on hereditary algebras. The main result is proved in section 5, while section 6
deals with the connection to cluster algebras, including necessary background.

The results of this paper have been presented at conferences in Uppsala
(June 2004), Mexico (August 2004) and Northeastern University (October 2004).

The first named author spent most of 2004 at the University of Leicester, and
would like to thank the Department of Mathematics, and especially Robert J. Marsh,
for their kind hospitality. We would like to thank the referee for pointing out an error
in an earlier version of this paper and Bernhard Keller and Otto Kerner for helpful comments
and conversations.

\section{Preliminaries}\label{prelim}

\subsection{Finite-dimensional algebras}

In this section let $\L$ be a finite dimensional $K$-algebra, where $K$ is a field.
Then $1_{\L} = e_1 + e_2 + \cdots + e_n$, where all $e_i$ are primitive idempotents.
We always assume that $\L$ is basic, that is, $\L e_i \not \simeq \L e_j$ when $i \neq j$.
There are then (up to isomorphism) $n$ indecomposable projective $\L$-modules, given by $\L e_i$, and
$n$ simple modules, given by $\L e_i / \r e_i$, where $\r$ is the Jacobson radical of $\L$.

If $K$ is an algebraically closed field, then 
there is a finite quiver $Q$, such that $\L$ is isomorphic to $KQ/I$, where
$KQ$ is the path-algebra, and $I$ is an admissible ideal, that is there is some $m$, such that
$\r^m \subseteq I \subseteq \r^2$.
We call $Q$ the quiver of $\L$. In case $\L$ is hereditary, the ideal $I$ is $0$.

The category $\mod \L$ of finite dimensional left $\L$-modules is an abelian category having 
almost split sequences. In case $\L$ is hereditary
there is a translation functor $\tau$, which is defined
on all modules with no projective (non-zero) direct summands.
Here $D$ denotes the ordinary duality for finite-dimensional algebras.

The bounded derived category of $\L$, denoted $D^b(\mod \L)$, is a triangulated category, with
suspension given by the shift-functor $[1]$, which is an autoequivalence.
We denote its inverse by $[-1]$.
In this paper, we only consider derived categories of hereditary algebras $H$. They  have
an especially nice structure, since the indecomposable objects are
given by shifts of indecomposable modules. In this case we also have a translation functor 
$\tau \colon  D^b(\mod H) \to D^b(\mod H)$, extending the functor mentioned above.
We have almost split triangles $A \to B \to C \to $ in $D^b(\mod H)$,
where $\tau C =A$, for each indecomposable $C$ in $D^b(\mod H)$.
We also have the formula $\Hom_{\D}(X, \tau Y) \simeq D\Ext^1_{\D}(Y,X)$, see \cite{hap}.
Let $H$ be a hereditary finite-dimensional algebra. Then a module $T$ in $\mod H$ is called
a {\em tilting module} if $\Ext^1_H(T,T) = 0$ and $T$ has, up to isomorphism, $n$ indecomposable direct summands.
The endomorphism ring $\End_H(T)^{\op}$ is called a {\em tilted algebra}.

See \cite{ars} and \cite{r} for further information on the representation theory of finite dimensional algebras and
almost split sequences.

\subsection{Approximations}

Let $\E$ be an additive category, and $\X$ a full subcategory.
Let $E$ be an object in $\E$. If there is an object $X$ in $\X$, and a map $f \colon X \to E$, such that
for every object $X'$ in $\X$ and every map $g \colon X' \to E$, there is a map 
$h \colon X' \to X$, such that $g = fh$, then $f$ is called a right $\X$-approximation \cite{as}.
The approximation map $f \colon X \to E$ is called {\em minimal} if no non-zero direct summand
of $X$ is mapped to $0$.
The concept of (minimal) left $\X$-approximations is defined dually.
If there is a field $K$, such that $\Hom_{\E}(X,Y)$ is finite dimensional over $K$, for all $X,Y \in \E$,
and if $\X = \add M$ for an object $M$ in $\E$, 
then (minimal) left and right $\X$-approximations always exist.

\subsection{Cluster categories and cluster-tilted algebras}
We remind the reader of the basic definitions and results from \cite{bmrrt}.
Let $H$ be a hereditary algebra, and let $\D = D^b(\mod H)$ be the bounded derived category.

The \emph{cluster category} is defined as the orbit category $\C_H = \D / F$, where $F = \tau^{-1}[1]$.
The objects of $C_H$ are the same as the objects of $\D$, but maps are given by
$$\Hom_{\C}(X,Y) = \amalg_i \Hom_{\D}(X,F^i Y).$$
Let $Q \colon \D \to \C$ be the canonical functor.
We often denote $Q(X)$ by $\widehat{X}$, and use the same notation for maps.
Let $\widehat{X}$ be an indecomposable object in the cluster-category.
We call $\mod H \vee \add H[1] = \mod H \vee H[1]$ the \emph{standard domain}.
There is 
(up to isomorphism) a unique 
object $X$ in $\mod H \vee H[1] \subseteq \D$ such that $Q(X) = \widehat{X}$.

Assume $X_1, X_2$ are indecomposable in the standard domain, then
a map $\widehat{f}\colon \widehat{X_1} \to \widehat{X_2}$, 
can uniquely be written as a sum of maps $\widehat{f_1}+ \widehat{f_2} + \cdots + \widehat{f_r}$,
such that $f_i$ is in $\Hom_{\D}(X_1, F^{d_i}X_2)$, 
for integers $d_i$. 
In this case $d_i$ is called
the {\em degree} of $\widehat{f_i}$.

The following summarizes properties about cluster categories that will be freely
used later.

\begin{thm}
Let $H$ be a hereditary algebra, and $\C_H$ the cluster category of $H$.
Then 
\begin{itemize}
\item[(a)]{$\C_H$ is a Krull-Schmidt category and $Q$ preserves indecomposable objects;}
\item[(b)]{$\C_H$ is triangulated and $Q$ is exact;}
\item[(c)]{$\C_H$ has AR-triangles and $Q$ preserves AR-triangles.}
\end{itemize}
\end{thm}

\begin{proof}
(b) is due to Keller \cite{k}, while (a) and (c) are proved in
\cite{bmrrt}.
\end{proof}

Let us now fix a hereditary algebra $H$, and assume it has, up to isomorphism, $n$ simple 
modules.
A \emph{cluster tilting object} (or for short; tilting object) in the 
cluster category is an object $T$ with $\Ext^1_{\C}(T,T) = 0$,
and with $n$ non-isomorphic indecomposable direct summands. For an object $X$ in any additive
category, we let $\add X$ denote the smallest full additive subcategory closed under direct sums
and containing $X$. Then two tilting objects $T$ and $T'$ are said to be
equivalent if and only if $\add T = \add T'$. We only consider
tilting objects up to equivalence, and therefore we always assume that
if $T = \amalg_i T_i$ is a tilting object, with each $T_i$ indecomposable,
then $T_i \not \simeq T_j$ for $i \neq j$.

There is a natural embedding of the module category into the bounded derived category, which extends
to an embedding of the module category into $\C$.
It was shown in \cite{bmrrt} that the image of a tilting module in $\mod H$ is a tilting object in
$\C_H$. It was also shown that if we choose a tilting object $T$ in $\C_H$, then there is
a hereditary algebra $H'$ and an equivalence $D^b(H') \to D^b(H)$, such that
$T$ is the image of a tilting module, under the embedding of $\mod H'$ into $\C_{H'} \simeq C_H$.

If $\bar{T} \amalg X$ is a tilting object, and $X$ is indecomposable, then
$\bar{T}$ is called an almost complete tilting object.

The following was shown in \cite{bmrrt}.

\begin{thm}\label{triangles}
Let $\bar{T}$ be an almost complete tilting object in $\C_H$. Then there are
exactly two complements $M$ and $M^{\ast}$. There are uniquely defined non-split
triangles 
$$M^{\ast} \to B \to M \to ,$$
and
$$M \to B' \to M^{\ast} \to .$$
The maps $ B \to M$ and $B' \to M^{\ast}$ are minimal 
right $\add \bar{T}$-approximations, and
the maps $M^{\ast} \to B$ and $M \to B'$ are minimal left
$\add \bar{T}$-approximations.
\end{thm}

The endomorphism ring $\End_{\C}(T)^{\op}$ of a tilting object in $\C$ is called a \emph{cluster-tilted algebra}.
Using the notation of Theorem \ref{triangles}, we want to 
compare the quivers of the endomorphism rings 
$\G = \End_{\C}(\bar{T} \amalg M)^{\op}$ and $\G' = \End_{\C}(\bar{T} \amalg M^{\ast})^{\op}$.

\subsection{Matrix mutation}

Let $X =(x_{ij})$ be an $n \times n$-skew-symmetric matrix with integer entries. 
Choose $k \in \{1,2,\dots,n\}$ and define a new matrix $\mu_k(X) = X' = (x'_{ij})$ by 
$$x'_{ij} = \begin{cases} 
-x_{ij} & \text{if $k=i$ or $k=j$,} \\
x_{ij} + \frac{\abs{x_{ik}}x_{kj} + x_{ik} \abs{x_{kj}}}{2} & \text{otherwise.}
\end{cases}
$$
The matrix $\mu_k(X) = X'$ is called the mutation of $X$ in direction $k$, and one can show
that
\begin{itemize}
\item[-]{$\mu_k(X)$ is skew-symmetric, and}
\item[-]{$\mu_k(\mu_k(X))= X$.}
\end{itemize}

Matrix mutation appears in the definition of cluster algebras by Fomin and Zelevinsky \cite{fz1}.

\subsection{Main result}

At this point, we have the necessary notation to state the main result of this paper.
There are no loops in the quiver of a cluster-tilted algebra \cite{bmrrt}, and we
also later show that there are no (oriented) cycles of length two. 
It follows that
one can assign to $\G$ a skew-symmetric integer matrix $X_{\G}$.
Actually, there is a 1--1 correspondence between the skew-symmetric integer matrices and
quivers with no loops and no cycles of length two.
Fixing an ordering of the vertices of the quiver, this 1--1 correspondence
determines mutations $\mu_k$ also on finite quivers (with no loops and no cycles of length 2).
The following will be proved in Section 5.
The notation is as earlier in this section, especially $K$ is algebraically closed.

\begin{thm}
Let $\bar{T}$ be an almost complete tilting object with complements
$M$ and $M^{\ast}$ and let $\G = \End_{\C}(\bar{T} \amalg M)^{\op}$ and 
$\G' = \End_{\C}(\bar{T} \amalg M^{\ast})^{\op}$. Let $k$ be the
vertex of $\G$ corresponding to $M$. Then the quivers $Q_{\G}$ and
$Q_{\G'}$, or equivalently the matrices
$X_{\G} = (x_{ij})$ and $X_{\G'} = (x'_{ij})$, are related by the formulas
$$x'_{ij} = \begin{cases} 
-x_{ij} & \text{if $k=i$ or $k=j$,} \\
x_{ij} + \frac{\abs{x_{ik}}x_{kj} + x_{ik} \abs{x_{kj}}}{2} & \text{otherwise.}
\end{cases}
$$
\end{thm}

This is the central result from which the connections with cluster algebras mentioned in the 
introduction follow.

\section{Factors of cluster-tilted algebras}

In this section, our main result is that for any cluster-tilted algebra $\G$,
and any primitive idempotent $e$, the factor-algebra $\G/ \G e \G$ is in 
a natural way also a cluster-tilted algebra. This will give us a powerful 
reduction-technique, which is of independent interest, and which we
use in the proof of our main result in this paper.

Suppose that $\G$ is the endomorphism algebra of a tilting object $T$ in
the cluster category corresponding to a hereditary algebra $H$.
The main idea of the proof is to show that if we localise $\D^b(\mod H)$
at the smallest thick subcategory containing a fixed indecomposable
summand $M$ of $T$, then we obtain a category triangle-equivalent to
the derived category of
a hereditary algebra $H'$. The factor-algebra $\G/ \G e \G$
(where $e$ is the primitive idempotent of $\G$ corresponding to $M$) is then
shown to be isomorphic to the endomorphism algebra of a tilting object in the
cluster category corresponding to $H'$.


\subsection{Localisation of triangulated categories} \label{localisationtri}

We review the basics of localisation in triangulated categories.
Let $\T$ be a triangulated category.
A subcategory $\M$ of $\T$ is called a \emph{thick subcategory} of
$\T$ if it is a full triangulated subcategory of $\T$ closed under
taking direct summands.

When $\M$ is a thick subcategory of $\T$, one can form a new
triangulated category
$\T_{\M}=\T/\M$, and there is a canonical exact functor
$L_{\M} \colon \T \to \T_{\M}$.
See \cite{rickard} and \cite{verdier} for details.

For every $M'$ in $\M$, we have $L_{\M}(M') = 0$, and $L_{\M}$ is
universal with respect to this property.
We also have the following.

\begin{lem}\label{maps}
Assume $\T$ is a triangulated category, and $\M$ is a thick
subcategory of $\T$.
Then, for any map $f$ in $\T$ we have $L_{\M}(f) = 0$ if and only if
$f$ factors through an object in $\M$.
\end{lem}

We will need the following result of Verdier \cite[Ch. 2, 5-3]{verdier}:

\begin{prop}\label{verdieriso}
Let $\T$ be a triangulated category with thick subcategory $\M$, and
let $\T_{\M}$ be the quotient category with quotient functor
$L_{\M}:\T\to \T_{\M}$.
Fix an object $Y$ of $\T$.
Then every morphism from an object of $\M$ to $Y$ is zero if and only
if for every object $X$ of $\T$ the canonical map
$$\Hom_{\T}(X,Y)\to \Hom_{\T_{\M}}(L_{\M}(X),L_{\M}(Y))$$
is an isomorphism.
\end{prop}

In particular, we note that this implies that $L_{\M}$ is fully
faithful on the full subcategory of $\T$ with objects given by those
objects of $\T$ which have only zero morphisms from objects of $\M$.

\subsection{Equivalences of module categories} \label{BongartzHappel}

Let $H$ be a hereditary algebra and $M$ an indecomposable $H$-module
with $\Ext^1_H(M,M)= 0$. Then there is (up to isomorphism) a unique module $E$ 
with the following properties:

\begin{itemize}
\item[B1)]{$E$ is a complement of $M$ (that is, $E \amalg M$ is a tilting module)}.
\item[B2)]{For any module $X$ in $\mod H$, we have that $\Ext^1_H(M,X) = 0$ implies
also $\Ext^1_H(E,X) = 0$}.
\end{itemize}

This is due to Bongartz \cite{bon}, and the module $E$ is sometimes
known as the \emph{Bongartz-complement} of $M$. 
For a module $X$ in  $\mod H$, 
we denote by $X^{\perp}$ the full subcategory of $\mod H$ with objects $Y$ satisfying
$\Ext^1_H(X,Y)=0$. 
If $T$ is a tilting
module, then it is well-known
that $T^{\perp} = \Fac T$, where $\Fac T$ is the full subcategory of all modules that are factors of
objects in $\add T$.
Note that B2) can be reformulated as $M^{\perp} = {(M \amalg E)^{\perp}}$.

The following result can be found in \cite{hap} and \cite{hris}.

\begin{prop} \label{happelequivalence}
\begin{itemize}
\item[(a)]{Assume $M$ is an indecomposable non-projective $H$-module with
$\Ext^1_H(M,M)= 0$, and
let $E$ be the complement as above.
Then the endomorphism ring
$H' = \End_H(E)^{\op}$ is hereditary, and $\Hom_H(M,E) = 0$.} 
\item[(b)]{
Let $\U$ denote the full
subcategory of $\mod H$ with objects $X$ satisfying
$$\Hom_H(M,X)=0= \Ext^1_H(M,X).$$
Then $\U$ is an exact subcategory of $\mod H$ and the functor
$\Hom_H(E,-)$ from $\mod H$ to $\mod H'$ restricts to an exact
equivalence between $\U$ and $\mod H'$.}
\end{itemize}
\end{prop}

We note that the above result does not hold in general in the case
when $M$ is projective. For example, consider the quiver of
type $A_3$ with vertices $1,2$ and $3$ and
arrows from $1$ to $2$ and $2$ to $3$.
Let $M=P_2$. Then $E=P_1\oplus P_3$
and $\End_H(E)^{op}$ has three indecomposable objects while $\U$ has only
two. The only other complement of $M$ is $E'=P_1\oplus (P_2/P_3)$.
Then $\Hom_H(E',P_3)=0$ although $P_3$ lies in $\U$ and is non-zero.
So also in this case the functor $\Hom_H(E',\ )$ from $\U$ to $\mod \End_H(E')^{\op}$
is not an equivalence.
However, we will need the following result which is along similar lines
for the case when $M$ is projective.

\begin{lem}\label{projectiveequivalence}
Let $M$ be an indecomposable projective $H$-module with corresponding
idempotent $e_M \in H$. Let $H'=H/{He_M H}$. 
\begin{itemize}
\item[(a)]{We have $\Tor_1^H(H',U)=0$,
for any object $U$ in $\U$, where $\U$ is as defined above.}
\item[(b)]{We have that
$\U$ is an exact subcategory of $\mod H$ and the functor $H'\otimes_H -$ 
from $\mod H$ to $\mod H'$ restricts to an exact equivalence between $\U$
and $\mod H'$.}
\end{itemize}
\end{lem}

\begin{proof}
We have that $\U$ is an exact subcategory of $\mod H$ as in
Proposition~\ref{happelequivalence}. It is easy to see that the
functor $H'\otimes_H -$ is an equivalence between $\U$ and $\mod H'$.
To see that it is exact, we consider the following projective resolution
of $H'$ as a right $H$-module:
$$0\to He_M H\to H\to H'\to 0.$$
Applying $-\otimes_H U$ to this sequence, where $U$ is an object in
$\U$, we obtain (part of) the long exact sequence:
$$\Tor_1^H(H,U)\to \Tor_1^H(H',U)\to He_M H\otimes_H U\to
H\otimes_H U\to H'\otimes_H U\to 0.$$
Since $H$ is projective, $\Tor_1^H(H,U)=0$. We also have
$$He_M H\otimes_H U=H\otimes_H He_MU=0$$
since $e_M U=0$. It follows that $\Tor^H_1(H',U)=0$ and hence that
$H'\otimes_H-$ is an exact functor on $\U$.
\end{proof}

\subsection{Localising with respect to an exceptional module}\label{loc}
Fix a hereditary algebra $H$, and
an indecomposable module $M$ in $\mod H$, with $\Ext^1_H(M,M) = 0$.

\begin{lem}
Let $\M = \add \{ M[i] \}_{i \in \Z}$. Then $\M$ is a thick subcategory in $D^b(\mod H)$.
\end{lem}

\begin{proof}
Straightforward from the fact that any map between indecomposable
objects in $\M$ is either zero or an isomorphism.
\end{proof}

Let $\D = D^b(\mod H)$, let $\D_{\M}$ be the category
obtained from $\D$ by localising with respect to $\M$, and let
$L_{\M} \colon \D \to \D_{\M}$ be the localisation functor. 
Note that $\U$ is the full subcategory of $\mod H$ consisting
of modules $X$ with $\Hom_{\D}(M,X[i]) = 0$ for all $i$.

\begin{thm}\label{both}
Let $H$ be a hereditary algebra with $n$ simple modules up to isomorphism.
Let $M$ be an indecomposable $H$-module with
$\Ext_H^1(M,M) = 0$, and let $\M$ denote the thick subcategory
generated by $M$. Then $\D_{\M}$ is equivalent to the derived category of a 
hereditary algebra with $n-1$ simple modules (up to isomorphism).
\end{thm}

To prove this we show that $\D_{\M}$ is equivalent to the subcategory
$\D_0 = \add \{X[i] \in \D \mid X \in \U, i \in \Z \}$ of  $\D$, and that $\D_0$
is equivalent to the derived category of a 
hereditary algebra with $n-1$ simple modules.
This is the content of the following three propositions.
We usually denote the object $L_{\M}(X)$ by $\widetilde{X}$.

\begin{prop} \label{localequiv}
In the setting of Theorem \ref{both}, the localisation functor
$L_{\M}$ induces an equivalence $\D_0 \to \D_{\M}$.
\end{prop}

\begin{proof}
First note that by Proposition \ref{verdieriso} we have that $L_{\M} \colon 
\D_0 \to \D_{\M}$ is fully faithful.
Any object in $\D_{\M}$ is of the form $L_{\M}(X)$ for some object $X$ in
$\D$.
Let $\widetilde{X}$ be an arbitrary object in $\D_{\M}$ (where $X$ is in $\D$).
Then consider the minimal right $\M$-approximation $M_X \to X$, and
the induced triangle $M_X \to X \to X_0 \to$. 
It is clear that $\widetilde{X} = \widetilde{X_0}$.
We claim that $X_0$ is in $\D_0$, that is $\Hom_{\D}(M, X[i]) = 0$
for all $i$. To see this, consider the long exact sequence
obtained by applying $\Hom_{\D}(M,\ )$ to the triangle $M_X \to X \to X_0 \to $.
For any $i$, the
map $\Hom_{\D}(M,M_X[i]) \to \Hom_{\D}(M,X[i])$ is an epimorphism, since $M_X \to X$ is a
right $\M$-approximation. The map is injective since any element in $\Hom_{\D}(M,M_X[i])$
is either zero or an isomorphism.
Thus, $X_0$ is in $\D_0$.
This completes the proof that $L_{\M}$ induces an equivalence $\D_0 \to \D_{\M}$.
\end{proof}

The next result is an extension of Proposition \ref{happelequivalence} 
to the setting of derived categories.

\begin{prop}
In the setting of Theorem \ref{both}, assume $M$ is non-projective.
Let $E$ be the Bongartz-complement of $M$, and let $H' = \End_H(E)^{op}$.
Then $\RHom(E,\ )$ induces an equivalence $\D_0 \to \D' = D^b(\mod H')$.
\end{prop}

\begin{proof}
Recall that 
$\U \subset M^{\perp} = (M \amalg E)^{\perp}$.
This implies that for $X \in \U$, we have that $\RHom(E,X)$ is concentrated in
degree zero with zero-term $\Hom_H(E,X)$.
Since $\Hom_H(E,\ )$ is a dense functor from $\U$ to 
$\mod H'$, and $\RHom(E,\ )$ commutes with $[1]$, it follows that
$\RHom(E,\ )$ restricted to $\D_0$ is dense.

Assume $X,Y$ are indecomposable objects in the same degree in $\D_0$. By the above it 
now follows directly from Proposition \ref{happelequivalence} that 
$$\Hom_{\D}(X,Y) \simeq \Hom_{\D'}(\RHom(E,X), \RHom(E,Y)).$$

We also need to show that $\Hom_{\D}(X,Y[1]) \simeq \Hom_{\D'}(\RHom(E,X), \RHom(E,Y[1]))$.
For this note that by Proposition \ref{happelequivalence}, the equivalence $\Hom_H(E,\ ) \colon \U \to \mod H'$
is exact, and that the embedding $\U \hookrightarrow \mod H$ is exact.
This implies that 
\begin{align*}
\Hom_{\D}(X,Y[1]) \simeq & \Hom_{\U}(X,Y[1])  \\
\simeq & \Hom_{\D'}(\Hom_H(E,X), \Hom_H(E,Y)[1])  \\
\simeq & \Hom_{\D'}(\RHom(E,X),\RHom(E,Y)[1])  \\
\simeq & \Hom_{\D'}(\RHom(E,X),\RHom(E,Y[1])).
\end{align*}
Thus the restriction of $\RHom(E,\ )$ to $\D_0$ is fully faithful. 
This completes the proof. 
\end{proof}

\begin{prop}
In the setting of Theorem \ref{both}, assume $M$ is projective. Assume 
$M \simeq He_M$ for the primitive idempotent $e_M$ in $H$
and let $H' =H / He_M H$. 
Then $\LL(H' \otimes_H - )$ induces an equivalence $\D_0 \to \D' = D^b(\mod H')$.
\end{prop}

\begin{proof}
First recall from Lemma \ref{projectiveequivalence} that 
$\Tor_1^H(H',U) = 0$ for any $U$ in $\U$.
This means that the image $\LL(H' \otimes_H U)$ is just $H' \otimes_H U$
concentrated in degree $0$.

It now follows that $\LL(H' \otimes_H -)$ restricted to $\D_0$ is dense, by using
that the functor $H' \otimes_H - \colon \U \to \mod H'$ is dense
and that $\LL(H' \otimes_H - )$ commutes with $[1]$.

Assume $X,Y$ are indecomposable objects in the same degree in $\D_0$. It
follows from Lemma \ref{projectiveequivalence} that 
$$\Hom_{\D}(X,Y) \simeq \Hom_{\D'}(\LL(H' \otimes_H X), \LL(H'
\otimes_H Y)).$$

We need also to show that 
$$\Hom_{\D}(X,Y[1]) \simeq \Hom_{\D'}(\LL(H' \otimes_H X), \LL(H'
\otimes_H Y[1])).$$
For this recall that the embedding of $\U$ into $\mod H$ is exact, and that
$H' \otimes_H -$ is exact on $\U$ by Lemma \ref{projectiveequivalence}. 
Thus it follows that:

\begin{align*}
\Hom_{\D}(X,Y[1]) \simeq &\Hom_{\U}(X,Y[1])  \\
\simeq &\Hom_{\D'}(H' \otimes_H X, H' \otimes_H Y [1]) \\ 
\simeq &\Hom_{\D'}(\LL(H' \otimes_H X),\LL(H' \otimes_H Y)[1]) \\
\simeq &\Hom_{\D'}(\LL(H' \otimes_H X),\LL(H' \otimes_H Y[1])).
\end{align*}

This shows that the functor is fully faithful and finishes the proof.
\end{proof}

For the remainder of this section, we view the induced equivalence between
$\D_{\M}$ and $\D'$ as an identification.

\subsection{The factor construction}

As before, let $M$ be an indecomposable $H$-module with
$\Ext^1_H(M,M)=0$, where $H$ is hereditary, and let $E$ be the Bongartz
complement of $M$. We investigate the image of an arbitrary complement
$\overline{T}$ of $M$ under the functor $L_{\M}$.
For an object $X$ in $\D$, we use the notation $\widetilde{X} = L_{\M}(X)$,
as before. Note that $L_{\M}(\overline{T})=
L_{\M}(T)=\widetilde{T}$.

\begin{lem}\label{orto}
Let the notation be as above.
\begin{itemize}
\item[(a)]{$L_{\M}(\overline{T})= \widetilde{T}$ is in $\mod H'\vee H'[1]$.
}
\item[(b)]{
$\Hom_{\D'}(\widetilde{T},\widetilde{T}[1])=0$.}
\end{itemize}
\end{lem}

\begin{proof}
(a) 
Let $f \colon M'\to \overline{T}$ be a minimal right
$\M$-approximation, and
consider the induced triangle:
\begin{equation} \label{approxtri}
M'\overset{f}\to \overline{T}\overset{g}\to U_T\to
\end{equation}
in $\D$. Since 
$\Hom_{\D}(M[-1],\overline{T})=0$, we have that
$M'$ is in $\add M$. 
It is clear that $\widetilde{T} \simeq \widetilde{U_T}$.
Now, as in the proof of Proposition \ref{localequiv}, we get that
$U_T$ is in $\D_0$. Here it is clear that $U_T = U_1 \amalg U_2[1]$, where
$U_1 = \Coker f$ and $U_2 = \Ker f$ are in $\U$. 
It is clear that $\widetilde{U_1}$ and $\widetilde{U_2}$ are $H'$-modules. 
We only need to show that $\widetilde{U_2}$ is projective.
For an arbitrary $U$ in $\U$, we have that $\Ext^1_H(U_2,U) = 0$, since $\Ext^1_H(M,U) = 0$
and $U_2$ is a submodule of $M$.
Using that $\U$ is an exact subcategory of $\mod H$, and that the equivalence $\U \to \mod H'$
is also exact, it follows that $\widetilde{U_2}$ is projective in $\mod H'$.
Hence $\widetilde{T} \simeq \widetilde{U_T}$ is in $\mod H'\vee H'[1]$.  

(b) Using again the triangle (\ref{approxtri}) we obtain the long exact sequence
$$\Hom_{\D}(\overline{T},\overline{T}[1])\to \Hom_{\D}(\overline{T},U_T[1])\to
\Hom_{\D}(\overline{T},M'[2]).$$
Hence, $\Hom_{\D}(\overline{T},U_T[1])= 0$. Now, by Proposition \ref{verdieriso},
it follows that $\Hom_{\D'}(\widetilde{T}, \widetilde{U_T[1]})= 0$ since $U_T[1]$ is
in $\D_0$, \sloppy and hence
$\Hom_{\D'}(\widetilde{T}, \widetilde{T}[1])= 0$.
%
%
\end{proof}


Denote as before by $F$ the functor $\tau^{-1}[1] \colon \D \to \D$. When
it is not clear which derived category $\D$ we are dealing with, we will denote this functor
by $F_{\D}$ and the functor $\tau^{-1}$ by $\tau_{\D}^{-1}$.

\begin{lem}\label{orto2} 
Let $H$ be a hereditary algebra, and let $X$ be an object
in $\D$ such that 
$X$ is in $\mod H \vee H[1]$. Then $\Hom_{\D}(X,X[1]) = 0$ if and only if
$\Ext^1_{\C_H}(X,X)=0$.
\end{lem}

\begin{proof}
Assume $X$ is in $\mod H \vee H[1]$, and let $\widehat{X}$ be the image of $X$ in the cluster
category $\C_H$ of $H$. Then $\Ext^1_{\C}(\widehat{X}, \widehat{X}) \simeq \Hom_{\D}(X,X[1]) 
\amalg D\Hom_{\D}(X,X[1])$. This follows from
$\Hom_{\D}(X, F^{-1}X[1]) = \Hom_{\D}(X, \tau X) \simeq D\Hom_{\D}(X,X[1])$ and 
the easily checked fact that $\Hom_{\D}(X, F^{i}X[1]) =  0$, whenever $i \not \in \{-1,0 \}$.
\end{proof}

Combining these lemmas, and using that a tilting $H$-module induces a
tilting object in the cluster category \cite[3.3]{bmrrt}, we obtain
the following.

\begin{prop}
Let $T = M \amalg \bar{T}$ be a tilting $H$-module as before.
Then the image $\widehat{T}$ of $\widetilde{T}$ in the cluster
category $\C_{H'}$ is a tilting object.
\end{prop}

\begin{proof}
Since $T$ is a tilting $H$-module, the triangulated category generated by $T$
is $\D$. Hence the triangulated category generated by $\widetilde{T}$ is $\D'$.
By Lemmas \ref{orto} and \ref{orto2} we have $\Ext^1_{\C_{H'}}(\widehat{T}, \widehat{T})=0$.
By the proof of \cite[Thm. 3.3]{bmrrt}, we have that $\widetilde{T}$ can be viewed as a 
direct summand of a tilting $H''$-module, for
some hereditary algebra $H''$ derived equivalent to $H'$. We then have
$\Ext^1_{H''}(\widetilde{T}, \widetilde{T})=0$, and hence
$\Hom_{\D'}(\widetilde{T}, \widetilde{T}[i]) = 0$ for all $i \neq 0$. Since
$\D'$ is the triangulated category generated by $\widetilde{T}$, it follows that
$\widetilde{T}$ is a tilting complex by definition, and consequently a tilting $H''$-module.
Then $\widehat{T}$ is a tilting object in $\C_{H'}$.
\end{proof}

We can now complete the main result of this section.
Let $e$ be the idempotent in $\G$, such that $\G e \simeq \Hom_{\C}(T,M)$.

\begin{thm}\label{factor}
With the above notation, there is a natural isomorphism 
$\G/ \G e \G \simeq \End_{\C_{H'}}(\widehat{T})^{\op}$.
\end{thm}

The remainder of this section will be devoted to proving this theorem.
Since the cluster category is defined using the functor
$F=\tau^{-1}[1]$, we need to compare $\widetilde{\tau_{\D}^{-1}(X)}$
and $\tau^{-1}_{\D'} \widetilde{X}$ for an indecomposable object $X$
in $\D$. In general
$\widetilde{\tau^{-1}_{\D} X} \not \simeq \tau^{-1}_{\D'} \widetilde{X}$,
but with extra conditions on $X$, sufficient for our purposes,
everything behaves nicely. The next result has been generalised by Keller, with a 
simpler proof \cite{k}, but for completeness we include our proof here.

\begin{lem}\label{commutes}
Let $X$ be an indecomposable object in $\D_0\subset \D$.
Then $\widetilde{X}$ is indecomposable and
$\widetilde{\tau^{-1}_{\D} X}  \simeq \tau^{-1}_{\D'} \widetilde{X}$.
\end{lem}

\begin{proof}
Since $L_{\M} \colon \D_0 \to \D_{\M}$ is an equivalence, we have that $\widetilde{X} = L_{\M}(X)$
is indecomposable. Let $f \colon X \to C$ be a minimal left almost
split map in $\D$. We want to show
that $\widetilde{f} \colon \widetilde{X} \to \widetilde{C}$ has the same property. We first
show that $\widetilde{f} \colon \widetilde{X} \to \widetilde{C}$
is not a split monomorphism. Assume to the contrary that there is a map $g' \colon \widetilde{C}
\to \widetilde{X}$ in $\D'$ such that $g' \widetilde{f} = \id_{\widetilde{X}}$. Since $X$ is in $\D_0$,
there is a map $g \colon C \to X$ in $\D$ such that $\widetilde{g} =
g'$
by Proposition \ref{verdieriso}.
Then $gf \colon X \to X$
is an isomorphism since $\widetilde{g} \widetilde{f} \colon \widetilde{X} \to \widetilde{X}$ is an isomorphism
and $X$ is in $\D_0$. Hence $f \colon X \to C$ is a split
monomorphism, and we have a contradiction. Hence
$\widetilde{f}:\widetilde{X}\to\widetilde{C}$ is not a split monomorphism.

We next show that $\widetilde{f} \colon \widetilde{X} \to \widetilde{C}$ is left almost split.
Let $\widetilde{h} \colon \widetilde{X} \to \widetilde{Y}$ be a map in $\D'$ which is not an isomorphism,
with $Y$ indecomposable in $\D_0$ and hence $\widetilde{Y}$ indecomposable in $\D'$.
Let $h \colon X \to Y$ be the map in $\D$, inducing the map
$\widetilde{h} \colon \widetilde{X} \to \widetilde{Y}$. This map is unique by Lemma \ref{verdieriso}.
Since $f \colon X \to C$ is left almost split, there is some $s \colon C \to Y$ such that 
$sf =h$. Hence we have $\widetilde{s}\widetilde{f} = \widetilde{h}$ in $\D'$, showing that 
$\widetilde{f} \colon \widetilde{X} \to \widetilde{C}$ is left almost split.

We also want to show that $\widetilde{f} \colon \widetilde{X} \to \widetilde{C}$ is a left minimal map.
Assume to the contrary that there is an indecomposable direct summand $C_1$ of $C$, such
that when $f_1 \colon X \to C_1$ is the map induced by $f \colon X \to C$,
then $\widetilde{f_1} \colon \widetilde{X} \to \widetilde{C_1}$ is 0, but $\widetilde{C_1} \neq 0$.
Then we have a commutative diagram 
$$
\xy
\xymatrix{
X \ar[dr]_{u_1} \ar[rr]^{f_1} & & C_1 \\
& M' \ar[ur]_{v_1} &
}
\endxy
$$
with $M'$ in $\M$. Since $f_1 \colon X \to C_1$ is irreducible, we have that 
either $u_1 \colon X \to M'$ is a split monomorphism, so that $X$ is in $\M$, and hence
$\widetilde{X} = 0$, or $v_1 \colon M' \to C_1$ is a split epimorphism
and hence $C_1$ is in $\M'$, so that $\widetilde{C_1} = 0$. In both cases we have a contradiction.
Hence we can conclude that $\widetilde{f} \colon \widetilde{X} \to \widetilde{C}$ is a left minimal map.

Let $X \overset{f_1}{\to} C \to \tau^{-1} X\to$ be an almost split triangle in $\D$.
Then the induced triangle $\widetilde{X} \overset{\widetilde{f_1}}{\to} \widetilde{C} \to 
\widetilde{\tau^{-1} X}$ is almost split since $\widetilde{f} \colon \widetilde{X} \to \widetilde{C}$ is a minimal 
left almost split map. Hence we get $\tau^{-1}_{\D'} \widetilde{X} \simeq \widetilde{\tau^{-1}_{\D}X}$.
\end{proof}

Let $T_x$ be an indecomposable direct summand in $T$, 
not isomorphic to $M$. Let $M_x \to T_x$ be a 
minimal right $\add M$-approximation, and consider as before
the induced triangle $$M_x \to T_x \to U_x \to$$
in $\D$, where we know that $U_x$ is in $\D_0$
by the proof of Proposition \ref{localequiv}.
Thus, by applying the above lemma to each of the indecomposable direct summands of $U_x$, 
we obtain $\widetilde{\tau^{-1}_{\D} U_x} \simeq \tau^{-1}_{\D'} \widetilde{U_x}$,
and thus $\widetilde{F_{\D} U_x} \simeq F_{\D'} \widetilde{U_x}$.
It is also clear that $\widetilde{U_x} \simeq \widetilde{T_x}$.

Now, pick two (not necessarily different) 
indecomposable direct summands $T_a$ and $T_b$
of $\bar{T}$.
Construct the triangle $$M_b \to T_b \to U_b \to,$$
as above, and apply $F$ to it, to obtain the triangle
$$FM_b \to FT_b \to FU_b \to.$$
Apply $\Hom_{\D}(T_a,\ )$ to this triangle, to obtain
the long exact sequence
\begin{multline}\label{unmod}
\Hom_{\D}(T_a, FM_b) \to \Hom_{\D}(T_a, FT_b) \to  \\
\Hom_{\D}(T_a, FU_b) \to
\Hom_{\D}(T_a, FM_b[1]).
\end{multline}
The last term vanishes, since $T_a$ and $M_b$ are modules. 
Since $M_b \to T_b$ is a minimal right $\add M$-approximation,
it follows that $FM_b \to FT_b$ is a minimal right $\add FM$-approximation.
We have that $\Hom_{\D}(T_a, FU_b) \simeq \Hom_{\D}(T_a, FT_b)/(FM)$, where
for an object $Z$ we use the notation $\Hom(X,Y)/(Z)$ to denote the $\Hom$-space
modulo maps factoring through an object in $\add Z$.

We claim there is an exact sequence 
$$
\Hom_{\D}(T_a, FM_b)/(M) \to \Hom_{\D}(T_a, FT_b)/(M) \to 
\Hom_{\D}(T_a, FU_b)/(M) \to 0
$$
induced from
the exact sequence (\ref{unmod}). For this
it is sufficient to show that the kernel of the second map
is contained in the image of the first.
So let $\alpha \in \Hom_{\D}(T_a, FT_b)/(M)$, and assume
there is a commutative diagram
$$
\xy
\xymatrix{
T_a \ar[r]^{\alpha} \ar[dr]^{\beta_1} & FT_b \ar[r] & FU_b \\
& M' \ar[ur] &
}
\endxy
$$
for some $M'$ in $\add M$.
Since $\Hom_{\D}(M,FM_b[1]) = 0$, there is $\beta_2 \colon M' \to FT_b$, such that
$M' \overset{\beta_2}{\to} FT_b \to FU_b = M' \to FU_b$. In $\Hom_{\D}(T_a,FT_b)/(M)$ we have
$\alpha = \alpha - \beta_2 \beta_1$. By using the long exact sequence (\ref{unmod}), we obtain 
that $\alpha = \alpha - \beta_2 \beta_1$ 
factors through $FM_b \to FT_b$, so
the sequence is exact.
It follows from this that 
$\Hom_{\D}(T_a, FT_b)/(M \amalg FM) \simeq \Hom_{\D}(T_a, FU_b)/(M)$.

Let $f \colon M_1 \to FU_b$ be a minimal right $\M$-approximation, and complete to
a triangle $M_1 \overset{f}{\to} FU_b \to (FU_b)'$.
Applying $\Hom_{\D}(T_a, \ )$,
we get an exact sequence 
$$\Hom_{\D}(T_a, M_1 ) \to \Hom_{\D}( T_a, FU_b) \to \Hom_{\D}(T_a, (FU_b)' ) \to \Hom_{\D}(T_a,M_1[1]).$$ 
Since $U_b$ is in degree $0$ or $1$, then $FU_b$ is in degree $1$,$2$
or $3$, so $M_1$ is in degree $0$,$1$,$2$ or $3$. 
Hence the indecomposable direct summands of $M_1[1]$ are in degree at least 1, so that 
$\Hom_{\D}(T_a,M_1[1]) = 0$. Since a map $h \colon T_a \to FU_b$ factors through an object in $\M$
if and only if it factors through the minimal right $\M$-approximation of $FU_b$, we get the
isomorphism
$$\Hom_{\D}(T_a, F_{\D}U_b)/ (M) \simeq \Hom_{\D}(T_a, (F_{\D}U_b)').$$
We get that this is isomorphic to $\Hom_{\D'}(\widetilde{T_a}, \widetilde{F_{\D}U_b})$, since 
$(F_{\D}U_b)'$ is in $\D_0$. By Lemma \ref{commutes} this is isomorphic to $\Hom_{\D'}(\widetilde{T_a}, F_{\D'}\widetilde{T_b})$.
We thus obtain that
$$\Hom_{\D}(T_a,F_{\D}T_b)/{(M\amalg FM)}\simeq \Hom_{\D'}(\widetilde{T_a},F_{\D'}\widetilde{T_b}).$$

We have 
$\Hom_{\D}(T_a, T_b) / (M \amalg FM) \simeq \Hom_{\D}(T_a, T_b)/ (M)$. Consider again
the triangle $M_b  \overset{f_b}{\to} T_b \to U_b$ in $\D$, where
$f_b \colon M_b \to T_b$ is a minimal right $\M$-approximation. Applying
$\Hom_{\D}(T_a, \ )$ gives an exact sequence
$$\Hom_{\D}(T_a, M_b) \to \Hom_{\D}( T_a, T_b) \to \Hom_{\D}(T_a, U_b ) \to \Hom_{\D}(T_a,M_b[1]).$$ 
Since $M_b$ is a module, we have $\Hom_{\D}(T_a, M_b[1]) = 0$, and hence
\linebreak $\Hom_{\D}(T_a, T_b)/(M) \simeq \Hom_{\D}(T_a, U_b)$, which is
isomorphic to $\Hom_{\D'}(\widetilde{T_a}, \widetilde{U_b})$ by
Proposition \ref{verdieriso}. We obtain that:
$$\Hom_{\D}(T_a,T_b)/{(M\amalg FM)}\simeq \Hom_{\D'}(\widetilde{T_a},\widetilde{T_b}).$$

Therefore $\G / \G e \G = \Hom_{\D}(T,T) \amalg \Hom_{\D}(T,FT) /
(M \amalg FM) \simeq \Hom_{\D'}(\widetilde{T},\widetilde{T}) \amalg
\Hom_{\D'}(\widetilde{T}, F_{\D'} \widetilde{T})$ as vector spaces.
It is straightforward to check that the map is also a ring map.
Theorem \ref{factor} is proved.

\subsection{Comparison with tilted algebras}

We give an example showing that a result similar to Theorem \ref{factor} does
not hold for tilted algebras. We would like to thank Dieter Happel for 
providing us with this example. 
There is a tilting module for the path algebra of a Dynkin quiver of type $D_5$, such
that the corresponding tilted algebra $\L$ has the quiver
$$
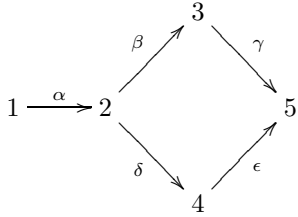
\xy
\xymatrix{
& & 3 \ar[dr]^{\gamma} & \\
1 \ar[r]^{\alpha} & 2 \ar[ur]^{\beta} \ar[dr]_{\delta} & & 5 \\
& & 4 \ar[ur]_{\epsilon} &
}
\endxy
$$
with relations $\alpha \beta = \beta \gamma - \delta \epsilon =0$.
If we let $e_4$ be the primitive idempotent corresponding to vertex $4$, then it
easy to see that $\L/\L e_4 \L$ is not tilted, since it has global dimension three.

It is well-known that the endomorphism-ring of a partial tilting
module is a tilted algebra. However, a similar result does
not hold for cluster-tilted algebras. An example of this is the path algebra 
of an oriented 4-cycle, modulo the cube of its radical. This 
is a cluster-tilted algebra of type $D_4$.

\section{Cluster-tilted algebras}\label{cluster-tilted}

In this section we apply the main result of the previous section to
show that (oriented) 
cycles in the quiver of a cluster-tilted algebra have length at least three.
Let $T_1 \amalg T_2 \amalg \dots \amalg T_n$ be a tilting object in
the cluster category $\C$. We denote by $\delta_k(T)$ the 
tilting object $T'$ obtained by exchanging $T_i$ with the 
second complement of $T_1 \amalg \cdots \amalg T_{i-1} \amalg T_{i+1} \amalg \cdots \amalg T_n$.
Let $\G = \End_{\C}(T)^{\op}$ and  $\G' = \End_{\C}(T')^{\op}$
be the corresponding cluster-tilted algebras. Passing from $\G$ to
$\G'$ depends on the choice of tilting object $T$. But we still write
$\overline{\delta}_k(\G) = \G'$, when either it is clear from the context
which tilting object $T$ gives rise to $\G$, or when this it is irrelevant.
We also say that $\G'$ is obtained from $\G$ by mutation at $k$.

From \cite{bmrrt} we know that all tilting objects in $\C_H$
can be obtained from performing a finite number of operations 
$\delta_k$ to $H$, where $H$ is the hereditary algebra considered as 
a tilting object in $\C_H$.

If $k$ is a source or a sink in the quiver of a hereditary algebra, then
mutation at $k$ coincides with so-called APR-tilting \cite{apr} (see \cite{bmr}), 
and the quiver of
the mutated algebra $\overline{\delta}_k(H)$ is obtained by reversing all 
arrows ending or starting in $k$.

\begin{lem}\label{ranktwo}
The cluster-tilted algebras of rank at most 2 are hereditary.
\end{lem}

\begin{proof}
This follows from the fact that any cluster-tilted algebra can be obtained 
by starting with a hereditary algebra, and performing a finite number
of mutations. If we start with a hereditary algebra $H$ of rank at most 2, the 
algebra obtained by mutating at one of the vertices is isomorphic to $H$. 
\end{proof}

\begin{prop}\label{noshortcycles}
The quiver of a cluster-tilted algebra has no 
loops and no cycles of length 2.
\end{prop}

\begin{proof}
This follows directly from combining 
Lemma \ref{ranktwo} with Theorem \ref{factor}.
\end{proof}

This was first proven by Gordana Todorov in case of finite representation type.

\section{Cluster-tilted algebras of rank 3}\label{rankthree}

In this section we specialize to connected hereditary algebras of rank 3, and 
the cluster-tilted algebras obtained from them.
We describe the possible quivers, and give some 
information on the relation-spaces. Later, this will be used to show our main result
for algebras of rank 3. In the proof of our main theorem, we use Theorem \ref{factor} 
to reduce to the case of rank 3. 
For hereditary algebras of finite representation type, there
is up to derived equivalence only one connected algebra of rank 3,
and thus up to equivalence only one cluster category $\C$.
In this case the technically involved results of this section reduce
to just checking {\em one case}: The only cluster-tilted algebra
of rank 3 which is not hereditary is given by a quiver which is a cycle of
length 3, and with the relations that the composition of any two arrows is zero.

\subsection{The quivers}

We consider quivers of the form
$$
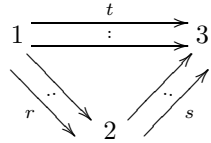
\xy
\xymatrix{1 \ar@<-0.5ex>[dr]_{.} \ar@<-2.5ex>[dr]^{.}_{r} \ar@<1ex>[rr]_{.}^{t} \ar@<-1ex>[rr]^{.} &  &3 \\
& 2 \ar@<-0.5ex>[ur]_{.} \ar@<-2.5ex>[ur]^{.}_{s} & 
}
\endxy
$$
where $r > 0$, $s > 0$ and $t \geq 0$
denote the number of arrows as indicated in the above figure. For short,
we denote such a quiver by $Q_{rst}$.

Up to derived equivalence, all connected finite dimensional hereditary algebras of rank 3 have
a quiver given as above. 
We first prove that factors of path-algebras of such quivers by admissible ideals are never 
cluster-tilted. The following is useful for this.

\begin{lem}\label{almost}
Let $H$ be a hereditary algebra of rank 3, and
assume $X, Y$ are indecomposable $H$-modules. Assume that
there exists an irreducible map $X \to Y$, and that 
$X \amalg Y$ is an almost complete tilting module. 
Then $X$ and $Y$ are either both preprojective or
both preinjective, and
there
exists a complement $Z$ such that $T = X \amalg Y \amalg Z$
has a hereditary endomorphism ring $\End_{\C}(T)^{\op} \simeq \End_H(T)^{\op}$.
\end{lem}

\begin{proof}
The first claim follows from the fact that
regular modules with no self-extensions are quasi-simple, since
$H$ has rank 3, and there does not exist an irreducible map between two quasi-simples.
Now, it is well-known that $X \amalg Y$ can be completed to a so-called {\em slice} \cite{r}, 
and the endomorphism ring $\End_H(T)^{\op}$ is hereditary.
Also, for a slice, $\Hom_{\D}(T,\tau^2 T)= 0$, so $\Hom_{\D}(T,FT) = \Hom_{\D}(T, \tau^{-1}T[1]) \simeq 
D\Hom_{\D}(T, \tau^2 T) = 0$. That is, $\End_{\C}(T)^{\op} \simeq \End_H(T)^{\op}$.
\end{proof}

\begin{lem}\label{nofactors}
If $\G$ is a
cluster-tilted algebra with quiver of type $Q_{rst}$, then $\G$ is hereditary.
\end{lem}

\begin{proof}
Assume $\G = \End_{\C}(T)^{\op}$, where $T = T_1 \amalg T_2 \amalg T_3$ is a tilting
module in $\mod H$.
For $H$ of finite type this is
easily checked, since there is, up to triangle-equivalence, only one 
cluster-category.
Therefore we can assume that $H$, and thus $\G$, is not of finite type.
Thus, we can choose $T$ and $H$ such that $T$ does not have both projective and injective 
direct summands. And by, if necessary, repeatedly applying $\tau^{-1}$ or $\tau$ we can assume that 
$T$, $\tau T$ and $\tau^{-1}T$ have no projective or injective direct summands.
There is a unique sink in $Q_{rst}$, and we assume that it corresponds to $T_3$.
Thus, we can assume that $\Hom_H(T_3, T_1 \amalg T_2) = 0$,
and $\Hom_{\D}(T_3, F(T))=0$. We use this to show that $\G$ is hereditary.

We now assume $\Hom_{\D}(T_3, \tau^{-1}T_i[1])=0$, and by the Auslander-Reiten formula
we have $\Hom_{H}(T_i, \tau^2 T_3)=0$, for $i = 1,2,3$.
We also have 
$\Hom_H(T_i, \tau T_3) = 0$, using the same formula.
Consider now the almost split sequence
\begin{equation}\label{almostsplit}
0 \to \tau T_3 \to X \to T_3 \to 0.
\end{equation}
We want to show that $X$ is in $\add T$.
Apply $\tau$ to (\ref{almostsplit}) to obtain
\begin{equation}\label{tauofit}
0 \to \tau^2 T_3 \to  \tau X \to \tau T_3 \to 0.
\end{equation}
By applying $\Hom_H(T,\ )$ to (\ref{tauofit}), and
using a long exact sequence argument, it follows that $\Hom_H(T, \tau X) = 0$, and so
$\Ext^1_H(X ,T)= 0$ by the AR-formula.

We next want to show that also 
$\Ext^1_H(T,X)=0$.
For $i = 1,2$ we show that
$\Ext^1_H(T_i,X) \simeq D \Hom_H(X, \tau T_i)=0$.
For this, apply $\Hom_H(\ ,\tau T_i)$ to the exact sequence (\ref{almostsplit})
to obtain the exact sequence
$$\Hom_H(T_3, \tau T_i) \to \Hom_H(X, \tau T_i) \to \Hom_H(\tau T_3, \tau T_i).$$
By the assumptions both
   \sloppy $\Hom_H(T_3, \tau T_i) \simeq D\Ext^1_H(T_i,T_3) = 0$ and
$\Hom_H(\tau T_3, \tau T_i) \simeq \Hom_H(T_3, T_i) =0$. Therefore also the middle term is zero.
By applying $\Hom_H(T_3,\ )$ to the exact sequence (\ref{almostsplit}), 
we obtain the long exact sequence
\begin{align*}
\Hom_H(T_3,\tau T_3) \to \Hom_H(T_3,X) \to \Hom_H(T_3,T_3) \to \\
\Ext^1_H(T_3,\tau T_3) \to \Ext^1_H(T_3,X) \to \Ext^1_H(T_3, T_3).
\end{align*}
Since (\ref{almostsplit}) is an almost split sequence, the 
map $\Hom_H(T_3,T_3) \to \Ext^1_H(T_3,\tau T_3)$
is an isomorphism. We have $ \Ext^1_H(T_3, T_3) = 0$, and
therefore
$\Ext^1_H(T_3,X) = 0$.
It follows from the same long exact sequence that $\Hom_H(T_3,X)= 0$, since
$\Hom_H(T_3,\tau T_3) \simeq D \Ext^1_H(T_3,T_3) =0$.

Now apply $\Hom_H(X,\  )$ to (\ref{almostsplit}), to obtain the long 
exact sequence 
$$\Ext^1_H(X,\tau T_3) \to \Ext^1_H(X,X) \to \Ext^1_H(X, T_3).$$
Since $\Ext^1_H(X,\tau T_3) \simeq D \Hom_H(T_3,X) = 0$, we get 
$\Ext^1_H(X,X) = 0$.
Thus $$\Ext^1_H(T \amalg X, T \amalg X)=0,$$ so $X$ is in $\add T$.
Since $T_1$ or $T_2$ must be a direct summand of $X$, there is an 
irreducible map in $\mod H$ from $T_1$ and/or $T_2$ to $T_3$.

If there are irreducible maps from both $T_1$ and $T_2$ to $T_3$, the claim
follows directly, since $T$  forms a slice in $\mod H$. 
Now assume there is an irreducible map $T_2 \to T_3$. 
For this, we need to discuss two cases.
If $t=0$, that is the quiver of $H$ is 
$$
\xy
\xymatrix{1 \ar@<1ex>[r]_{.}^{r} \ar@<-1ex>[r]^{.} &  2 \ar@<1ex>[r]_{.}^{s} \ar@<-1ex>[r]^{.} & 3 \
}
\endxy
$$
it can be easily seen that both complements to $T_2 \amalg T_3$
will form a slice, and this
will give a hereditary endomorphism ring.

Assume now $t>0$.
By Lemma \ref{almost}, the module $M \amalg T_2 \amalg T_3$ has a hereditary 
endomorphism ring, either for $M = T_1$ or for $M= T_1^{\ast}$.
Here $T_1^{\ast}$ is the second complement of $T_2 \amalg T_3$,
using the notation from Theorem \ref{triangles}.
If $M = T_1$, then we are done.
Assume therefore $M = T_1^{\ast}$.

In case $t>0$, there are two types of irreducible maps $X \to Y$: either
there exists an indecomposable module $Z$ with irreducible maps $X \to Z$ and $Z \to Y$, or no
such $Z$ exists. In the latter case, it is easily seen that 
both complements of $X \amalg Y$ give a hereditary endomorphism ring. 
In the first case, the AR-quiver looks like
$$
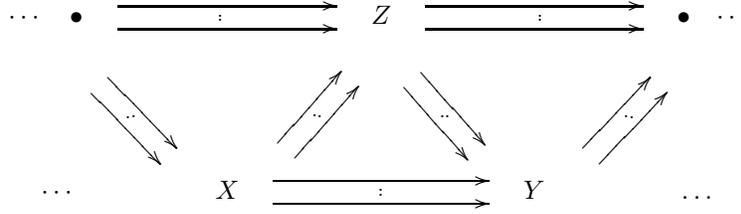
\xy
\xymatrix{
{\begin{matrix} & &  \\ \cdots & \bullet & \\ & & 
\end{matrix}}   \ar@<1ex>[rr]_{.} \ar@<-1ex>[rr]^{.} \ar@<-0.5ex>[dr]_{.} \ar@<-2.5ex>[dr]^{.}
& &   {\begin{matrix} & &  \\ & Z & \\ & & \end{matrix}} \ar@<1ex>[rr]_{.} \ar@<-1ex>[rr]^{.} 
\ar@<-0.5ex>[dr]_{.} \ar@<-2.5ex>[dr]^{.}
& &  {\begin{matrix} & &  \\ & \bullet & \cdots \\ & & \end{matrix}}  \\
{\begin{matrix} & &  \\ & \cdots & \\ & & \end{matrix}}  & 
{\begin{matrix} & &  \\ & X & \\ & & \end{matrix}}   \ar@<1ex>[rr]_{.} \ar@<-1ex>[rr]^{.} 
\ar@<-0.5ex>[ur]_{.} \ar@<-2.5ex>[ur]^{.} & & 
{\begin{matrix} & &  \\ & Y & \\ & & \end{matrix}}
\ar@<-0.5ex>[ur]_{.} \ar@<-2.5ex>[ur]^{.} & 
{\begin{matrix} & &  \\ & \dots & \\ & & \end{matrix}} 
}
\endxy
$$
One complement of $X \amalg Y$ is clearly $Z$, which gives a hereditary endomorphism ring.
The complement $Z^{\ast}$ is obtained from the triangle $Z \to B' \to Z^{\ast} \to$, where
$Z \to B$ is a minimal left $\add(X \amalg Y)$-approximation in $\C$. It is easy to see
that this approximation is just the left almost split map $Z \to \amalg Y$, where
$\amalg Y$ is a direct sum of copies of $Y$. This means that
there are non-zero maps $Y \to Z^{\ast}$ in $\C$. Consider the 
triangle $Z^{\ast} \to B \to Z \to$, then it is clear that $B \to Z$ is the right
almost split map $\amalg X \to Z$, where $\amalg X$ is a direct sum of copies of $X$,
and so there are non-zero maps $Z^{\ast} \to X$ in $\C$.
Thus there is a cycle in the quiver of $\End_{\C}(X \amalg Y \amalg  Z^{\ast})^{\op}$,
which gives a contradiction. 
This completes the proof for this case, and hence the proof of the lemma. 
\end{proof}

This has the following consequence.

\begin{cor}\label{cycles}
The quiver of a non-hereditary connected cluster-tilted algebra of rank 3 is of the form
$$
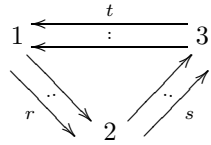
\xy
\xymatrix{1 \ar@<-0.5ex>[dr]_{.} \ar@<-2.5ex>[dr]^{.}_{r}  &  &
3 \ar@<1ex>[ll]_{.} \ar@<-1ex>[ll]^{.}_{t} \\
& 2 \ar@<-0.5ex>[ur]_{.} \ar@<-2.5ex>[ur]^{.}_{s} & 
}
\endxy
$$
with $r,s,t >0 $.
\end{cor}

\begin{proof}
Combine Lemma \ref{nofactors} with Propostion \ref{noshortcycles}.
\end{proof}

In view of this we refer to the cluster-tilted algebras of rank 3
which are non-hereditary as {\em cyclic cluster-tilted algebras}.

\subsection{The relations}

We first show that relations are homogeneous.

\begin{prop}\label{homogenous}
Let $\Gamma$ be a cluster-tilted algebra of rank 3 with Jacobson radical
$\underline{r}$. Then $\underline{r}^6 = 0$, and the relations are 
homogeneous.
\end{prop}

\begin{proof}
Without loss of generality we can assume that there is a tilting module
$T = X \amalg Y \amalg Z$ for a hereditary algebra $H$, such that
$\G = \End_{\C_H}(T)^{\op}$.

Using Corollary \ref{cycles} it is clear that we can assume
that the quiver of $\G$ has the form
$$
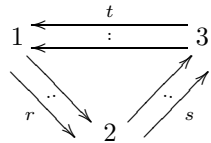
\xy
\xymatrix{1 \ar@<-0.5ex>[dr]_{.} \ar@<-2.5ex>[dr]^{.}_{r}  &  &
3 \ar@<1ex>[ll]_{.} \ar@<-1ex>[ll]^{.}_{t} \\
& 2 \ar@<-0.5ex>[ur]_{.} \ar@<-2.5ex>[ur]^{.}_{s} & 
}
\endxy
$$
with $r,s,t >0 $.

Let $\L  = \End_{H}(T)^{\op}$ be the corresponding tilted algebra.
There are no cycles in
the quiver of a tilted algebra. We can therefore assume that there is
a sink in the quiver of $\L$, and we assume that this vertex corresponds to
$Z$, that is, $\Hom_H(Z,X \amalg Y) = 0$. We assume $X,Y,Z$ correspond to
the vertices $1,2,3$, respectively. 
If $\widehat{h}$ is a non-zero map in $\Irr_{\add T}(Z,X)$, it must be of degree $1$,
that is, the lifting $h$ is in $\Hom_{\D}(Z,FX)$. Since this holds for all maps
in $\Irr_{\add T}(Z,X)$, any composition of $6$ arrows will
correspond to a map of degree $\geq 2$ from an indecomposable to itself,
and therefore must be the zero-map. This follows from the
fact that for any indecomposable module $M$, we have $\Hom_{\D}(M,F^2 M) = 0$.
This gives $\underline{r}^6 = 0$.

We can assume that at least one of the arrows (irreducible maps) $X \to Y$ and at 
least one of the arrows $Y \to Z$ are of degree 0.
Otherwise, the tilted algebra $\L$ would
not be connected.

Now let $\widehat{g}$ be a map in $\Irr_{\add T}(Y,Z)$. We want to show that 
it must be of degree 0.
Since $X \amalg Y$ is an almost complete tilting object in $\C_H$, there
are exactly two complements. Denote as usual the second one  by $Z^{\ast}$.
The complement $Z^{\ast}$ is either the image of a module or the image of
an object of the form $I[-1]$ for an injective indecomposable module $I$.
Furthermore, there is a triangle in $\C$
\begin{equation}\label{exchange}
Z^{\ast} \to Y^r \to Z \to ,
\end{equation}
for some $r \geq 0$, which can be lifted to 
a triangle 
$$F^i Z^{\ast} \overset{\big( \begin{smallmatrix} 
\alpha_1 \\ \alpha_2 \end{smallmatrix} \big) } \to 
Y^{r_1} \amalg (F^{-1}Y)^{r_2} \to Z \to $$
in $\D$ for some integer $i$ and with $r = r_1 + r_2$.
We need to show that $r_2 = 0$. It is sufficient to show that
the map $\alpha_2= 0$. 
We have $r_1 \neq 0$, and thus by minimality $\alpha_1 \neq 0$.
It is clear that if also $\alpha_2 \neq 0$,
then $i=0$ or $i= -1$.

Assume first $Z^{\ast} \simeq I[-1]$, then
$$\Hom_{\D}(I[-1],F^{-1}Y) = \Hom_{\D}(I, \tau Y) = 0,$$
so $i=0$ gives $\alpha_2 = 0$. On the other hand, it is
clear that $i = -1$ gives $\alpha_1 = 0$.

Assume now that $Z^{\ast}$ is the image of a module.
Then there is an exact sequence of modules
$$0 \to Z^{\ast} \to Y^s \to Z \to 0,$$
and since $\dim_k \Hom_{\C}(Z,Z^{\ast}[1]) =1$ (by \cite{bmrrt}), it follows that the triangle
(\ref{exchange})
is induced by this sequence, and thus $r_1 = s$ and $r_2 = 0$.

Now we show that all the irreducible maps $X \to Y$ in $\C_H$ 
are of degree 0. For this, consider
the almost complete tilting object $X \amalg Z$ in $C_H$,
with complements $Y$ and $Y^{\ast}$.
Consider the triangle 
$$Y^{\ast} \to X^t \to Y \to ,$$
and the preimage in $\D$,
$$F^i Y^{\ast}  \overset{\big( \begin{smallmatrix} 
\beta_1 \\ \beta_2 \end{smallmatrix} \big) } \to X^{t_1} \amalg (F^{-1}X)^{t_2} \to Y \to .$$

We need to show that $t_2 = 0$. 
The case where $Y^{\ast} \simeq I[-1]$ is completely similar as
for irreducible maps $Y \to Z$. In case $Y^{\ast}$ is the image of a module,
it is now more complicated since we have two possibilities. Either
there is an exact sequence in $\mod H$ of the form
$$0 \to Y^{\ast} \to X^u \to Y \to 0,$$ 
or there is an exact sequence of the form
$$0 \to Y \to Z^v \to Y^{\ast} \to 0.$$
If we are in the first case, we can use the same argument as for 
irreducible maps $Y \to Z$. If we are in the second case,
note that $\Hom_H(Y^{\ast},X) = 0$, since $\Hom_H(Z,X) = 0$. 
Thus, either $\beta_1= 0$ or $\beta_2= 0$
in our triangle. This completes the proof that
all irreducible maps $X \to Y$ are induced by module maps, and thus are of degree 0.

Given that $\underline{r}^6 = 0$, the only possibility for 
a non-homogeneous relation must involve maps in $\r^2 \setminus \r^3$ and maps
in $\r^5$. But, by our description of irreducible maps, this is not possible, because it would
involve a relation between maps of different degrees.
\end{proof}

Fix a cyclic cluster-tilted algebra of rank 3, and fix a vertex $k$. Let $\alpha$ be an
arrow ending in $k$, and $\beta$ an arrow starting in $k$. If $\beta \alpha = 0$, as an
element of the algebra, for any choice of $\alpha$ and $\beta$, then we call
$k$ a {\em zero vertex}.

\begin{prop}\label{zero}
Let $\G$ be a cyclic cluster-tilted algebra, and fix a vertex $k$. Then $k$ is a 
zero-vertex if and only if $\overline{\delta}_k(\G)$ is hereditary.
\end{prop}

\begin{proof}
We assume the quiver of $\G$ is
$$
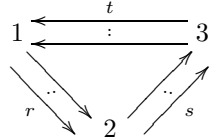
\xy
\xymatrix{1 \ar@<-0.5ex>[dr]_{.} \ar@<-2.5ex>[dr]^{.}_{r}  &  &
3 \ar@<1ex>[ll]_{.} \ar@<-1ex>[ll]^{.}_{t} \\
& 2 \ar@<-0.5ex>[ur]_{.} \ar@<-2.5ex>[ur]^{.}_{s} & 
}
\endxy
$$
Let $\G = \End_{\C}(T)^{\op}$, and let $T_i$ be the direct summand of $T$ corresponding to the vertex $i$.
Assume that $2$ is a zero-vertex.
Then it is clear that $\Hom_{\C}(T_1, T_3 ) = 0$, so the
quiver of $\overline{\delta}_2(\G)$ must be 
$$
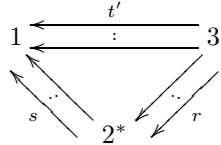
\xy
\xymatrix{1  &  & 3 \ar@<0.5ex>[dl]^{.} \ar@<2.5ex>[dl]_{.}^{r} \ar@<-1ex>[ll]^{.}_{t'}  \ar@<1ex>[ll]_{.} \\
& 2^{\ast} \ar@<0.5ex>[ul]^{.} \ar@<2.5ex>[ul]_{.}^{s} & 
}
\endxy
$$
with $t' \geq 0$.
Now $\overline{\delta}_2(\G)$ is hereditary, by Lemma \ref{nofactors}.

Conversely, assume $\overline{\delta}_k(\G)$ is hereditary. The quiver of $\overline{\delta}_k(\G)$ must be as above,
with $t' \geq 0$. This means $\Hom_{\C}(T_1, T_3 ) = 0$, so $2$ is a zero-vertex.
\end{proof}

\subsection{Kerner's Theorem}

The following result by Kerner \cite{ker} turns out to be crucial for the proof of 
the main theorem of this section. There is a more general version of this
theorem in \cite{ker}. We include a proof, for the convenience of the reader.
This proof is also due to Kerner, and we thank him for 
providing us with it.

\begin{thm}
Let $X,Y$ be regular indecomposable modules over a wild hereditary algebra $H$ of
rank 3. If $\Hom_H(X, \tau Y) = 0$, then also  $\Hom_H(X, \tau^{-1} Y) = 0$. 
\end{thm}

\begin{proof}
We first prove the following.
\begin{lem}\label{special}
Let $U$ be an indecomposable regular module over a wild hereditary algebra
of rank 3. Then $\Hom_H(U, \tau^2 U) \neq 0$.
\end{lem}

\begin{proof}
Assume first $\Ext^1_H(U,U) \neq 0$. By the AR-formula, then also $\Hom_H(U, \tau U) \neq 0$.
Assume now $\Hom_H(U, \tau^2 U) = 0$. Then also $\Ext^1_H(\tau U, U)= 0$ and,
by the Happel-Ringel lemma \cite{hr}, a non-zero map $f \colon U \to \tau U$ is either surjective
or injective. 
In either case, $g = \tau (f) \circ f \colon U \to \tau^2 U$ is non-zero.
This contradicts $\Hom_H(U, \tau^2 U) = 0$.

Now assume  $\Ext^1_H(U,U) = 0$. Then by \cite{Hos}, $U$ is quasi-simple. Thus,
there is an almost split sequence $0 \to \tau U \to V \to U \to 0$, where
$V$ is indecomposable, and by \cite{ker2} we have $\End_H(V) \simeq K$, while $\Ext^1_H(V,V) \neq 0$.
Applying $\Hom_H(U,\ )$ to the almost split sequence, we obtain the exact sequence
$$\Hom_H(U, \tau U) \to \Hom_H(U,V) \to \Hom_H(U, U) \to \Ext^1_H(U,\tau U)$$
Since $\Hom_H(U,U) \to \Ext^1_H(U,\tau U)$ is an isomorphism and $\Hom_H(U, \tau U) = 0$, we have that also 
$\Hom_H(U,V) = 0$. 
The long exact sequence obtained by 
applying $\Hom_H(\ , \tau U)$ to the almost split sequence, gives $\Hom_H(V, \tau U) = 0$.
Now, this gives $\Hom_H(V, \tau^2 U) \neq 0$, since there is an exact sequence
$$0 \to \Hom_H(V, \tau^2 U) \to \Hom_H(V, \tau V) \to \Hom_H(V, \tau U)$$
and the last term is zero. There is also the long exact sequence
$$0 \to \Hom_H(U, \tau^2 U) \to \Hom_H(V, \tau^2 U) \to \Hom_H(\tau U, \tau^2 U)$$
where the last term is zero. This proves $\Hom_H(U, \tau^2 U) \neq 0$.
\end{proof}

Let us now complete the proof of the theorem. Let $X,Y$ be regular indecomposable modules.
It suffices to show that $\Hom_H(X,Y) \neq 0$ implies $\Hom_H(X, \tau^2 Y) \neq 0$.
Let $z \colon X \to Y$ be a non-zero map. Then we can assume there is an indecomposable regular module
$U$, such that $z$ factors as $X \overset{p}{\to} U \overset{i}{\to} Y$, where
$p$ is surjective and $i$ is injective. Also $\tau^2 i \colon \tau^2 U \to \tau^2 Y$ is injective.
By Lemma \ref{special}, there is a non-zero map $f \colon U \to \tau^2 U$. The composition $\tau^2 i \circ f \circ p$
is non-zero. This completes the proof of the theorem.
\end{proof}

\subsection{The dimensions of relation-spaces}

Let $H$ be a connected hereditary algebra of rank 3. The following notation
is used for the rest of this section. Let $\bar{T}$ be an almost 
complete tilting object with complements $M$ and $M^{\ast}$, and assume there
are triangles as in Theorem \ref{triangles}.
Let $T = \bar{T} \amalg M$ and $T' = \bar{T} \amalg M^{\ast}$ 
and let $\G = \End_{\C}(T)^{\op}$ and  $\G' = \End_{\C}(T')^{\op}$. By now, we
know that the quiver of $\G$ is either 
$$
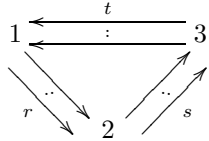
\xy
\xymatrix{1 \ar@<-0.5ex>[dr]_{.} \ar@<-2.5ex>[dr]^{.}_{r}  &  &
3 \ar@<1ex>[ll]_{.} \ar@<-1ex>[ll]^{.}_{t} \\
& 2 \ar@<-0.5ex>[ur]_{.} \ar@<-2.5ex>[ur]^{.}_{s} & 
}
\endxy
$$
with $r,s,t >0$ or the quiver $Q_{rst}$ 
$$
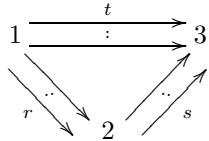
\xy
\xymatrix{1 \ar@<-0.5ex>[dr]_{.} \ar@<-2.5ex>[dr]^{.}_{r} \ar@<1ex>[rr]_{.}^{t} \ar@<-1ex>[rr]^{.} &  &3 \\
& 2 \ar@<-0.5ex>[ur]_{.} \ar@<-2.5ex>[ur]^{.}_{s} & 
}
\endxy
$$
with $r,s > 0$ and $t \geq 0$.
We let $M$ correspond to vertex $2$. Then $\bar{T} = T_{B} \amalg T_{B'}$ where $T_B$ corresponds to
the vertex $1$ and $T_{B'}$ to $3$. It is then clear that 
$B = (T_B)^r$ and $B' = (T_{B'})^s$. We label the vertices with the corresponding modules,
then the arrows represent irreducible maps in $\add T$.

We let $I$ denote the ideal such that $\G \simeq KQ/I$. 
In case $\G$ is cyclic,
we say that $\G$ is {\em balanced} at the vertex $2$ if
$$\dim((\Irr(T_B,M) \otimes \Irr(M,T_{B'}) \cap I) = t.$$
We will show that any vertex of a cyclic cluster-tilted algebra is either
balanced or a zero-vertex. 
We first discuss the algebras obtained by mutating hereditary algebras.

\begin{lem}\label{mutatehereditary}
Let $H$ be a hereditary algebra with quiver
$Q_{rst}$
where $r,s>0$ and $t \geq 0$. Then the following hold. 
\begin{itemize}
\item[(a)]{The cluster-tilted algebra $\G' = \overline{\delta}_2(H)$ is balanced at 
at the vertices $1$ and $3$.} 
\item[(b)]{The new vertex $2^{\ast}$ is a zero-vertex.}
\item[(c)]{The quiver of $\G'$ is
$$
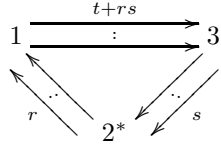
\xy
\xymatrix{1  \ar@<1ex>[rr]_{.}^{t + rs} \ar@<-1ex>[rr]^{.} &  & 3 
\ar@<0.5ex>[dl]^{.} \ar@<2.5ex>[dl]_{.}^{s} \\
& 2^{\ast} \ar@<0.5ex>[ul]^{.} \ar@<2.5ex>[ul]_{.}^{r} & 
}
\endxy
$$
}
\end{itemize}
\end{lem}

\begin{proof}
Part (b) and (c) follow directly from Lemma \ref{nofactors}
and Proposition \ref{zero}.
Let $P_i$ be the indecomposable projective $H$-module corresponding to vertex $i$, and $S_i$ 
the simple $H$-module $P_i/\r P_i$. Then $P_3 = S_3$ is simple. 
Consider $P_1 \amalg P_3$ as an almost complete tilting object.
There is an exact sequence $$0 \to P_2 \to (P_1)^s \to P_2^{\ast} \to 0,$$ such that
the induced triangle in $\C$
is the exchange-triangle of Theorem \ref{triangles}. Let $T' = P_1 \amalg P_2^{\ast} 
\amalg P_3$. Using the definition of $\tau$, one can show that
$S_2 = \tau P_2^{\ast}$, and thus $\Hom_{\D}(P_2^{\ast}, P_1)= 0$.
Since $2^{\ast}$ is a zero-vertex, $\Irr_{\add T'}(P_3, P_1) = \Hom_H(P_3, P_1)$,
with dimension $rs +t$.
We want to compute $\Irr_{\add T'}(P_2^{\ast}, P_3) \otimes \Irr_{\add T'}(P_3, P_1)
= \Hom_{\C}(P_2^{\ast},P_1) \simeq \Hom_{\D}(F^{-1}P_2^{\ast},P_1)$.
We have \sloppy $\Hom_{\D}(F^{-1}P_2^{\ast},P_1) =  \Hom_{\D}(\tau P_2^{\ast} [-1], P_1) = \Ext_H^1(S_2,P_1)$.
There is an exact sequence $$0 \to (P_3)^r \to P_2 \to S_2 \to 0.$$
Apply $\Hom_H(\ ,P_1)$ to it, to obtain the long exact sequence
$$0 \to \Hom_H(S_2, P_1) \to \Hom_H(P_2, P_1) \to \Hom_H((P_3)^r,P_1) \to \Ext^1_H(S_2,P_1) \to 0.$$
Since $\dim \Hom_H(P_2, P_1) = s$, and $\dim \Hom_H((P_3)^r,P_1) = (rs +t)r$, 
\sloppy we have $\dim (\Irr_{\add T'}(P_2^{\ast}, P_3) \otimes \Irr_{\add T'}(P_3, P_1)) = 
r(rs+t)-s$, thus $\dim (\Irr_{\add T'}(P_2^{\ast}, P_3) \otimes \Irr_{\add T'}(P_3, P_1) \cap I) =s$,
and $\G'$ is balanced at $3$.

Now apply $\Hom_H(P_3, \ )$ to the exact sequence
$ 0 \to P_2 \to (P_1)^s \to P_2^{\ast} \to 0$ to obtain the exact sequence
$$0 \to \Hom_H(P_3,P_2) \to \Hom_H(P_3,P_1^{s}) \to \Hom_H(P_3,P_2^{\ast}) \to 0.$$
Since $\dim  \Hom_H(P_3,P_2) = r$ and $\dim \Hom_H(P_3,P_1^{s}) = (rs +t)s$, we have
$$\dim (\Irr_{\add T'}(P_3,P_1) \otimes \Irr_{\add T'}(P_1, P_2^{\ast})) = \dim \Hom_{\C}(P_3,P_2^{\ast}) = 
(rs + t)s - r.$$ This means $\dim (\Irr_{\add T'}(P_3,P_1) \otimes \Irr_{\add T'}(P_1, P_2^{\ast}) \cap I) = r$,
and $\G'$ is balanced also at $1$.
\end{proof}

\begin{prop}\label{mutatecyclic}
Let $\G$ be a non-hereditary cluster-tilted algebra with quiver
$$
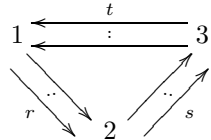
\xy
\xymatrix{1 \ar@<-0.5ex>[dr]_{.} \ar@<-2.5ex>[dr]^{.}_{r}  &  &
3 \ar@<1ex>[ll]_{.} \ar@<-1ex>[ll]^{.}_{t} \\
& 2 \ar@<-0.5ex>[ur]_{.} \ar@<-2.5ex>[ur]^{.}_{s} & 
}
\endxy
$$
\begin{itemize}
\item[(a)]{If $\G$ is balanced 
at the vertex $2$, then $\G' = \overline{\delta}_2(\G)$ is non-hereditary, and thus cyclic,
with quiver
$$
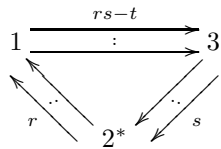
\xy
\xymatrix{1  \ar@<1ex>[rr]_{.}^{rs-t} \ar@<-1ex>[rr]^{.} &  & 3 
\ar@<0.5ex>[dl]^{.} \ar@<2.5ex>[dl]_{.}^{s} \\
& 2^{\ast} \ar@<0.5ex>[ul]^{.} \ar@<2.5ex>[ul]_{.}^{r} & 
}
\endxy
$$ 
It is balanced at the new vertex $2^{\ast}$. Each of the other vertices of $\G'$ is 
either balanced or a zero-vertex.}
\item[(b)]{If $\G$ has a zero-vertex at $2$, then $\overline{\delta}_2(\G)$ is hereditary with
quiver
$$
\xy
\xymatrix{1  &  & 3 \ar@<0.5ex>[dl]^{.} \ar@<2.5ex>[dl]_{.}^{r} \ar@<-1ex>[ll]^{.}_{t -rs}  \ar@<1ex>[ll]_{.} \\
& 2^{\ast} \ar@<0.5ex>[ul]^{.} \ar@<2.5ex>[ul]_{.}^{s} & 
}
\endxy
$$
}
\end{itemize}
\end{prop}

\begin{proof}
Part (b) follows from Propositions \ref{homogenous} and \ref{zero}.

To prove part (a), we adopt our earlier notation and conventions. Especially, 
$\G= \End_{\C}(T_B \amalg T_{B'} \amalg M)^{\op}$, and we have the quiver
$$
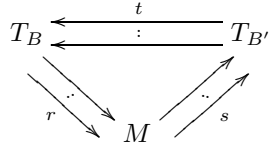
\xy
\xymatrix{T_B \ar@<-0.5ex>[dr]_{.} \ar@<-2.5ex>[dr]^{.}_{r}  &  &
T_{B'} \ar@<1ex>[ll]_{.} \ar@<-1ex>[ll]^{.}_{t} \\
& M \ar@<-0.5ex>[ur]_{.} \ar@<-2.5ex>[ur]^{.}_{s} & 
}
\endxy
$$
The quiver of the 
mutated algebra $\overline{\delta}_2(\G) = \End_{\C}(T_B \amalg T_{B'} \amalg M^{\ast})^{\op}$
is
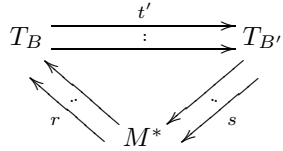
$$
\xy
\xymatrix{T_B  \ar@<1ex>[rr]_{.}^{t'} \ar@<-1ex>[rr]^{.} &  & T_{B'} 
\ar@<0.5ex>[dl]^{.} \ar@<2.5ex>[dl]_{.}^{s} \\
& M^{\ast} \ar@<0.5ex>[ul]^{.} \ar@<2.5ex>[ul]_{.}^{r} & 
}
\endxy
$$ 
Using that $\G$ is balanced at $2$, and Proposition \ref{homogenous},
it follows that $t' = rs -t$. Also by assumption, $M$ does not correspond to a zero-vertex, so there is 
at least one non-zero composition $T_B \to M \to T_{B'}$.
Therefore $rs -t > 0$.

We have $\dim (\Irr_{\add T'}(T_{B'},M^{\ast}) \otimes \Irr_{\add T'}(M^{\ast}, T_{B})) = \dim
\Irr_{\add T}(T_{B'},T_B) + \dim (\Irr_{\add T}(T_{B'},M) \otimes \Irr_{\add T}(M, T_{B}))= t+ 0$.
Therefore $\dim (\Irr_{\add T'}(T_{B'},M^{\ast}) \otimes \Irr_{\add T'}(M^{\ast}, T_{B}) \cap I) = rs-t$,
so $\G'$ is balanced at $2^{\ast}$.

We now proceed to show that for each of the vertices 1 and 3, $\G'$ is either 
balanced, or a zero-vertex.
We assume $T_{B'}$ is not a zero-vertex in $\G'$.

The tilted algebra $\L = \End_H(T_B \amalg T_{B'} \amalg M)^{\op}$ has a unique sink.
There is an induced total ordering on the triple $T_B, T_{B'}, M$, where the last 
element in the ordering
corresponds to the sink. Also, by considering the preimage of $M^{\ast}$ in the standard
domain of $\D$, the ordering can be
extended to the quadruple $B'$, $M^{\ast}$, $B$, $M$. Note that we get the following
four possible orderings 
\begin{itemize}
\item[-]{$(M, T_{B'}, M^{\ast}, T_B)$}
\item[-]{$(T_{B'}, M^{\ast}, T_B, M)$}
\item[-]{$(M^{\ast}, T_B, M, T_{B'})$}
\item[-]{$(T_B, M, T_{B'}, M^{\ast})$.}
\end{itemize}

First we show the claim for the vertex corresponding to $T_{B'}$.

\begin{lem}\label{magic}
Assume that $T_{B'}$ does not correspond to a zero-vertex in $\G'$ and
that $M^{\ast}$ is before $T_B$ in the above ordering. Then 
$\Hom_{\D}(T_B, \tau^{-1}M^{\ast})=0$. 
\end{lem}

\begin{proof}
Since $\overline{\delta}_2(\G)$ is not hereditary, we have $\Hom_H(M^{\ast},T_B) \neq 0$. 
Assume now that $\Hom_H(T_B,\tau^{-1}M^{\ast}) \neq 0$. Assume first that 
$T_B$ is regular, then $M^{\ast}$
is also regular. In case $H$ is tame, then 
there are at most two exceptional modules which are regular. This follows
from the fact that $H$ has three simples. 
But in case there are two exceptional modules which are regular, there
is an extension between them. 
This gives a contradiction.
In case $H$ is wild we can apply 
Kerner's Theorem, which says that $\Hom_H(T_B, \tau M^{\ast}) \neq 0$.
We have a contradiction, since $\Hom_H(T_B, \tau M^{\ast}) \simeq D\Ext^1_H(M^{\ast},T_B)= 0$.

If $B$ is a preprojective or a preinjective module, then $\Hom_H(M^{\ast},T_B) \neq 0$ and 
$\Hom_H(T_B, \tau^{-1}M^{\ast}) \neq 0$ implies that the map $M^{\ast} \to T_B$ is
irreducible in the module-category. Thus $M^{\ast} \amalg T_B$ can be complemented
to a tilting module with hereditary endomorphism ring. We have seen that the 
mutated algebra $\overline{\delta}_2(\G)$ is by
assumption not hereditary. This means that $T_{B'}$ must correspond to 
a zero-vertex, so we have a contradiction to $\Hom_H(B,\tau^{-1}M^{\ast}) \neq 0$
also for $T_B$ being preprojective or preinjective. 
\end{proof}

Now, let $M \to B' \to M^{\ast} \to $ be the usual triangle.
We recall that $\bar{T}= T_B \amalg T_{B'}$.
Let $\widetilde{\Hom}(T_B,M) = \Irr_{\add T}(T_B,M)$, let $\widetilde{\Hom}(T_B,B') = \Irr_{\add \bar{T}}(T_B,B')$
and let $\widetilde{\Hom}(T_B,M^{\ast}) = \Irr_{\add T'}(T_B,M^{\ast})$.
Then we claim that there is an exact sequence 
\begin{equation}\label{tildes}
0 \to \widetilde{\Hom}(T_B,M) \overset{\alpha}\to \widetilde{\Hom}(T_B,B') \to 
\widetilde{\Hom}(T_B,M^{\ast}) \to 0.
\end{equation}
It is clear from Proposition \ref{homogenous}, and the fact that
$$\Hom_{\D}(T_B,M) \to \Hom_{\D}(T_B,B') \to \Hom_{\D}(T_B,M^{\ast}) \to 0$$
is exact, that we only need to show that the map $\alpha$ is a monomorphism.
We first assume $M^{\ast}$ is a module.
To prove the claim for this case, we consider the four orderings 
on the quadruple $\{M,B',M^{\ast},B\}$.
For each case we show that a map in $\Irr_{\add T}(T_B,M)$ cannot factor via $M^{\ast}[-1]$ in $\C$. 
\vspace{2mm}
\\
\noindent 
$(M,T_{B'},M^{\ast},T_B)\colon$ In this case, a map in $\Irr_{\add T}(T_B,M)$ is of degree 
$1$. Assume the lifting is $f \colon \tau T_B[-1] \to M$.
There is a non-split exact sequence $0 \to M \to B' \to M^{\ast} \to 0$.
We have $\dim \Hom_{\C}(M^{\ast},M[1])= 1$, by \cite{bmrrt}, and therefore $\Hom_{\D}(FM^{\ast},M[1]) = 
\Hom_{\D}(\tau^{-1}M^{\ast},M) = 0$. Therefore, if $f\colon T_B \to M$ factors through $M^{\ast}[-1]$ in $\C$,
there must be a map $g \colon \tau T_B[-1] \to M$ in $\D$, such that there is a commutative diagram
$$
\xy
\xymatrix{
& \tau T_B[-1] \ar[d]_{f} \ar[dl]_{g} \\
M^{\ast}[-1] \ar[r] & M
}
\endxy
$$
By Lemma 
\ref{magic}, we have that $\Hom_{\D}(\tau T_B,M^{\ast})= \Hom_{\D}(T_B, \tau^{-1}M^{\ast})= 0$, 
and thus we obtain $f=0$.
\vspace{2mm}
\\
\noindent
$(T_{B'}, M^{\ast}, T_B, M)$ or $(M^{\ast}, T_B, M, T_{B'}) \colon$
In these cases, a map in $\Irr_{\add T}(T_B,M)$ 
is of degree $0$. Assume the lifting of it is $f \colon T_B \to M$.
The preimage of $M^{\ast}$ in $\D$ is a module in these cases, so
a factorization of $f$ must be of the form
$$
\xy
\xymatrix{
& T_B \ar[d]_{f} \ar[dl]_{g} \\
\tau^{-1}M^{\ast} \ar[r] & M
}
\endxy
$$
Lemma \ref{magic} gives $f=0$.
\vspace{2mm}
\\
\noindent
$(T_B, M, T_{B'}, M^{\ast})\colon$ 
In this case the preimage of $M^{\ast}$ in $\D$ is either a module or
$P[1]$, for an indecomposable projective $H$-module $P$.
In both cases, a map in $\Irr_{\add T}(T_B,M)$ 
is a map of degree $0$. Assume the lifting is $f\colon T_B \to M$.
In the first case there is a non-split exact sequence
$0 \to M \to B' \to M^{\ast} \to 0$. Therefore, 
since $\dim \Hom_{\C}(M^{\ast}, M[1]) = 1$, we have $\Hom_{\D}(\tau^{-1}M^{\ast},M) = 0$.
Since $\Hom_{\D}(T_B, M^{\ast}[-1])=0$, we must have $f=0$, if $f$
factors as below
$$
\xy
\xymatrix{
& T_B \ar[d]_{f} \ar[dl]_{g} \\
M^{\ast}[-1] \ar[r] & M.
}
\endxy
$$
Assume $M^{\ast} \simeq P[1]$, with $P$ projective. 
If $\Hom_{\D}(T_B,P) \neq 0$, then $T_B$ is also projective. Therefore
$\Hom_{\D}(T_B, P[1]) = 0$, and thus $\Hom_{\C}(T_B, M^{\ast})= 0$, which
means that $T_{B'}$ is a zero-vertex in $\G'$, a contradiction.
Then $\Hom_{\D}(T_B,P) = 0$, but this means that $f \colon  T_B \to M$ factors
through $M^{\ast}[-1] \simeq P$ only for $f =0$.
Thus, the map $\alpha$ is a monomorphism, and the sequence (\ref{tildes}) is exact.

Thus, $\dim \widetilde{\Hom}(T_B,M^{\ast}) = \dim \widetilde{\Hom}(T_B,B')- \dim \widetilde{\Hom}(T_B,M)= 
t' s - r$. 
This means that $\dim (\Irr_{\add T'}(T_B,T_{B'}) \otimes \Irr_{\add T'}(T_{B'}, M^{\ast}) \cap I ) = r$,
so $\G'$ is balanced at the vertex $3$, corresponding to $T_{B'}$.

We now show that $\G'$ is balanced at the vertex $1$, corresponding to $T_B$,
or $1$ is a zero-vertex. Assume it is not a zero-vertex.
We have the dual version of Lemma \ref{magic}.

\begin{lem}
Assume that $T_{B}$ does not correspond to a zero-vertex and
that $M^{\ast}$ is before $T_B$ in the above ordering. Then 
$\Hom_{\D}(M^{\ast}, \tau T_{B'})=0$. 
\end{lem}

\begin{proof}
Similar to the proof of Lemma \ref{magic}.
\end{proof}

Now, consider the triangle $$M^{\ast} \to B \to M \to.$$ We need to show that
there is an exact sequence 
\begin{equation}\label{tildes2}
0 \to \widetilde{\Hom}(M,T_{B'}) \to \widetilde{\Hom}(B,T_{B'}) \to 
\widetilde{\Hom}(M^{\ast},T_{B'}) \to 0,
\end{equation}
where $\widetilde{\Hom}(M,T_{B'})= \Irr_{\add T}(M,T_{B'})$, while $\widetilde{\Hom}(B,T_{B'})= 
\Irr_{\add \bar{T}}(B,T_{B'})$ and $\widetilde{\Hom}(M^{\ast},T_{B'}) = \Irr_{\add T'}(M^{\ast},T_{B'})$.
The proof of this is parallel to the proof for the sequence (\ref{tildes}), and therefore omitted.
Using the exact sequence (\ref{tildes2}), one obtains
that $\G'$ is balanced at the vertex $1$.
\end{proof}

We summarize the results of this section.
\begin{thm}\label{sum}
Let $\G$ be a cluster-tilted algebra of rank 3.
\begin{itemize}
\item[(a)]{$\G$ is either
hereditary, or it is cyclic.}
\item[(b)]{
If $\G$ is cyclic, then each vertex of $\G$ is either balanced or a zero-vertex.}
\item[(c)]{Let $k$ be a vertex of $\G$, let $\overline{\delta}_k(\G)$ be the mutation in direction $k$,
and let $k^{\ast}$ be the new vertex of $\overline{\delta}_k(\G)$.
Then there are the following possible cases:
\begin{itemize}
\item[I.]{Both $\G$ and $\overline{\delta}_k(\G)$ are hereditary.}
\item[II.]{$\G$ is hereditary, while $\overline{\delta}_k(\G)$ is cyclic with a zero-vertex at $k^{\ast}$,}
\item[III.]{$\G$ is cyclic with a zero-vertex at $k$, and $\overline{\delta}_k(\G)$ is hereditary, or}
\item[IV.]{$\G$ is cyclic and balanced at $k$, and $\overline{\delta}_k(\G)$ is cyclic and balanced at $k^{\ast}$}.
\end{itemize}
}
\end{itemize}
\end{thm}

\begin{proof}
This follows directly from the previous results in this section, and the fact that
all cluster-tilted algebras can be obtained by starting with a hereditary algebra,
and then performing a finite number of mutations \cite{bmrrt}, \cite{bmr}.
\end{proof}

The above Theorem is very easily verified for algebras of finite type,
as indicated in the introduction of this section.

\section{Mutation}

As mentioned in the introduction, in view of Proposition \ref{noshortcycles} 
it is possible to assign to 
a cluster-tilted algebra $\G$ a skew-symmetric matrix $X_{\G} = (x_{ij})$.
More precisely, if there is at least one arrow from $i$ to $j$
in the quiver of the endomorphism-algebra $\G$, let
$x_{ij}$ be the number of arrows from $i$ to $j$. If there are no arrows between 
$i$ and $j$, let $x_{ij} = 0$. Let $x_{ij} = - x_{ji}$ otherwise.

Now let $\bar{T}$ be an almost complete tilting object with complements
$M$ and $M^{\ast}$. Let $T = \bar{T} \amalg M$, let $T' = \bar{T} \amalg M^{\ast}$,
let $\G = \End_{\C}(T)^{\op}$ and
let $\G' = \End_{\C}(T')^{\op}$. Then we want to show that
the quivers of $\G$ and $\G'$ are related by the cluster-mutation formula. 
We use the results of Section \ref{rankthree} to show this for
cluster-tilted algebras of rank 3, and Theorem \ref{factor} to extend to the
general case.

\begin{thm}\label{mutate} 
Let $H$ be a hereditary algebra, and let 
$\bar{T},M , M^{\ast}, \G$ and $\G'$ be as above.
Then the quivers of $\G$ and $\G'$, or equivalently the matrices 
$X_{\G}$ and $X_{\G'}$, are related by the 
cluster mutation formula.
\end{thm}

\begin{proof}
First, assume $H$ has rank 3. In case $H$ is not connected, the  claim is easily checked.
Assume $H$ is connected. 
Fix $k$, the vertex where we mutate. By Theorem \ref{triangles},
it is clear that $x_{ik}' = - x_{ik}$ for $i= 1,2,3$, and that
$x_{kj}' = - x_{kj}$ for $j= 1,2,3$.

Now assume $i \neq k$ and $j \neq k$. 
By Theorem \ref{sum}, there are four possible cases.
\vspace{2mm}
\\
\noindent 
Case I: This happens if and only if $k$ is a source or a sink. In this case it is
clear that either $x_{ik}= 0$ or $x_{kj}= 0$. 
For $i \neq k$ and $j \neq k$,
it is clear that $x_{ij} = x_{ij}'$, since in this case mutation at $k$ is the
same as so-called APR-tilting at $k$. Thus the formula holds.
\vspace{2mm}
\\
\noindent
Case II: Since $k$ is now not a source or a sink, we can assume $\G$ is the path algebra of
$$
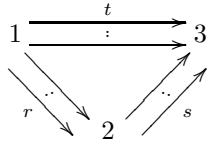
\xy
\xymatrix{1 \ar@<-0.5ex>[dr]_{.} \ar@<-2.5ex>[dr]^{.}_{r} \ar@<1ex>[rr]_{.}^{t} \ar@<-1ex>[rr]^{.} &  &3 \\
& 2 \ar@<-0.5ex>[ur]_{.} \ar@<-2.5ex>[ur]^{.}_{s} & 
}
\endxy
$$
where $r>0$ and $s > 0$ and $t \geq 0$ and with $k=2$. Then, by Lemma \ref{mutatehereditary}, the quiver
of $\overline{\delta}_2(\G)$ is 
$$
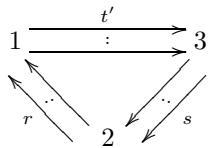
\xy
\xymatrix{1  \ar@<1ex>[rr]_{.}^{t'} \ar@<-1ex>[rr]^{.} &  & 3 
\ar@<0.5ex>[dl]^{.} \ar@<2.5ex>[dl]_{.}^{s} \\
& 2 \ar@<0.5ex>[ul]^{.} \ar@<2.5ex>[ul]_{.}^{r} & 
}
\endxy
$$
with $t' = r s + t$.
So $x_{13}' = t' = rs +t$, and 
$$x_{13} + \frac{\abs{x_{12}}x_{23} + x_{12} \abs{x_{23}}}{2} = t + \frac{rs + rs}{2} = t + rs. $$

\noindent
Case III:  We assume that the quiver of $\G$ is
$$
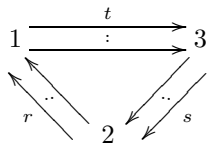
\xy
\xymatrix{1  \ar@<1ex>[rr]_{.}^{t} \ar@<-1ex>[rr]^{.} &  & 3 
\ar@<0.5ex>[dl]^{.} \ar@<2.5ex>[dl]_{.}^{s} \\
& 2 \ar@<0.5ex>[ul]^{.} \ar@<2.5ex>[ul]_{.}^{r} & 
}
\endxy
$$
By Proposition \ref{mutatecyclic}, the quiver of $\G'$ is 
$$
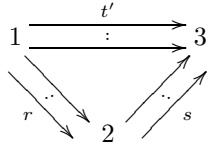
\xy
\xymatrix{1 \ar@<-0.5ex>[dr]_{.} \ar@<-2.5ex>[dr]^{.}_{r} \ar@<1ex>[rr]_{.}^{t'} \ar@<-1ex>[rr]^{.} &  &3 \\
& 2 \ar@<-0.5ex>[ur]_{.} \ar@<-2.5ex>[ur]^{.}_{s} & 
}
\endxy
$$
with $t' = t-rs$.
That is $x_{13}' = t-rs$, and $$x_{13} + \frac{\abs{x_{12}}x_{23} + x_{12} \abs{x_{23}}}{2} =
t + \frac{\abs{-r} (-s) + (-r) \abs{-s}}{2} = t-rs,$$ 
and the formula holds.

\noindent 
Case IV: We assume the quiver of $\G$ is the same as in case III. Now the quiver of
$\G'$ is
$$
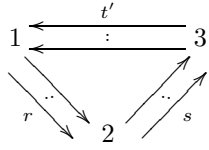
\xy
\xymatrix{1 \ar@<-0.5ex>[dr]_{.} \ar@<-2.5ex>[dr]^{.}_{r}  &  &
3 \ar@<1ex>[ll]_{.} \ar@<-1ex>[ll]^{.}_{t'} \\
& 2 \ar@<-0.5ex>[ur]_{.} \ar@<-2.5ex>[ur]^{.}_{s} & 
}
\endxy
$$
where $t' = rs-t$. That is $x_{13}'= -t' = t-rs$, while $$x_{13} + \frac{\abs{x_{12}}x_{23} 
+ x_{12} \abs{x_{23}}}{2} =
t + \frac{\abs{-r} (-s) + (-r) \abs{-s}}{2} = t-rs,$$
thus the formula holds true also in this case.

Now, assume that $H$ has arbitrary rank. 
Fix $k$, the vertex where we mutate. By Theorem \ref{triangles},
it is clear that $x_{ik}' = - x_{ik}$ for any value of $i$, and that
$x_{kj}' = - x_{kj}$ for any value of $j$.
Assume now that $k \neq i$ and $k \neq j$.
Let $e_i, e_j, e_k$ be the primitive idempotents in $\G$ corresponding to the vertices 
$i,j,k$ of the quiver of $\G$. Assume $1_{\G}= f +e_i +e_j +e_k$.
and let $\G_{\red} = \G / \G f \G$. Let $e_i, e_j, e_k^{\ast}$ be the primitive 
idempotents corresponding to the vertices $i,j,k^{\ast}$ of the quiver of $\G'$.
Assume $1_{\G'}= f' +e_i +e_j +e_{k^{\ast}}$.
and let $\G'_{\red} = \G' / \G' f' \G'$.
It is clear that the number of arrows from $i$ to $j$ in the quiver of 
$\G_{\red}$ is $x_{ij}$ and the number of arrows from $i$ to $j$ in the quiver of 
$\G'_{\red}$ is $x_{ij}'$. So, by the first part of the proof, $x_{ij}$ and
$x_{ij}'$ are related by the matrix mutation formula.
\end{proof}

\section{Connections to cluster algebras}

Our main motivation for studying matrix mutation for quivers/matrices associated
with tilting objects in cluster categories is the connection to cluster algebras.
In this section we explain how Theorem \ref{mutate} gives such a connection.
In order to formulate our result we first need to give a short introduction to a special type of 
cluster algebras \cite{fz1}, relevant to our setting \cite{bfz}.
See also \cite{fz2} for an overview of the theory of cluster algebras.

Let $\F=\Q(u_1, \dots, u_n)$ be the field of rational functions in indeterminates
$u_1, \dots , u_n$, let $\underline{x} = \{x_1, \dots, x_n\} \subset F$ be a transcendence
basis over $\Q$, and $B= (b_{ij})$ an $n \times n$ skew-symmetric integer matrix.
A pair $(\underline{x},B)$ is called a {\em seed}.
The {\em cluster algebra} associated to the seed $(\underline{x},B)$ is by definition
a certain subring $\A(\underline{x},B)$ of $\F$, as we shall describe.
Given such a seed $(\underline{x},B)$ and some $i$, with $1 \leq i \leq n$,
define a new element of $x_i'$ of $\F$ by
$$x_i x_i' = \prod_{j; b_{ji}>0} x^{b_{ji}}+ \prod_{j; b_{ji}<0}  x^{-b_{ji}}.$$ 
We say that $x_i, x_i'$ 
form an {\em exchange pair}. We obtain a new transcendence basis 
$\underline{x'} = \{x_1, \dots, x_n\} \cup \{x_i'\} \setminus \{x_i \}$ of $\F$.
Then define a new matrix $B' = (b_{ij}')$ associated with $B$ by
$$b'_{ij} = \begin{cases} 
-b_{ij} & \text{if $k=i$ or $k=j$,} \\
b_{ij} + \frac{\abs{b_{ik}}b_{kj} + b_{ik} \abs{b_{kj}}}{2} & \text{otherwise.}
\end{cases}
$$
The pair $(\underline{x'},B')$ is called the {\em mutation} of the seed $(\underline{x},B)$
in direction $i$, written $\mu_i(\underline{x},B) = (\underline{x}',B')$.
Let $\S$ be the set of seeds obtained by iterated mutations of $(\underline{x},B)$
(in all possible directions). The set of {\em cluster variables} is by definition the
union of all transcendence bases appearing in all the seeds in $\S$, and
the cluster algebra $\A(\underline{x},B)$ is the subring of 
$\F$ generated by the cluster variables. The transcendence bases appearing in the seeds are 
called {\em clusters}.

As mentioned earlier,
there is a 1--1 correspondence between
finite quivers with no loops and no oriented cycle of length two and
skew-symmetric integer matrices (up to reordering the columns). The vertices of the quiver 
of a matrix $B=(b_{ij})$ are $1, \dots, n$, and there are $b_{ij}$ arrows from $i$ to $j$
if $b_{ij}> 0$.
The cluster algebra is said to be {\em acyclic} if there is some seed where
the quiver associated with the matrix has no oriented cycles \cite{bfz}. We take the
corresponding seed as an initial seed.
In this case, let $H= KQ$ be the hereditary path algebra associated with an initial seed 
$(\underline{x},B)$. Let $\C = \C_H$ be the corresponding cluster category, and let $T$
be a tilting object in $\C$.
Similar to the above we can associate with $T$ a {\em tilting seed} $(T, Q_T)$, where
$Q_T$ is the quiver of the endomorphism algebra $\End_{\C}(T)^{\op}$. Let
$T_1, \dots, T_n$ be the non-isomorphic indecomposable direct summands of $T$.
Fix $i$, and let as before $\delta_i(T) = T'$ be the tilting object of
$\C$ obtained by exchanging $T_i$ with $T_i^{\ast}$ (using our earlier notation
from Theorem \ref{triangles}).
Define mutation of $(T, Q_T)$ in direction $i$ to be given by $\delta_i(T,Q_T) = (T', Q_{T'})$.

We now want to associate tilting seeds with seeds for acyclic cluster algebras. We first associate
$(H[1], Q_H)$ with a fixed initial seed $(\underline{x},B)$, where $Q$ is the quiver for 
$B$ and $H= KQ$. Let $ (\underline{x}',B')$ be some seed. We then have
$ (\underline{x}',B') = \mu_{i_t}\cdots \mu_{i_1} (\underline{x},B)$ for some ordered sequence
$(i_1, \dots i_t)$. There are in general several such sequences, and we choose one of minimal length.
Associated with $(\underline{x},B)$ is the sequence of length $0$, that is the empty set
$\emptyset$.
We define $\alpha((\underline{x},B), \emptyset)= (H[1],Q_H)$, and
$\alpha((\underline{x}',B'), (i_1, \dots,i_t))= \delta_{i_t} \cdots \delta_{i_1} (H[1], Q_H) =
(T', Q_{T'})$.
Fix an ordering on the cluster variables in the cluster $\underline{x}= \{x_1, \dots x_n\}$
of the chosen initial seed and choose a corresponding indexing for the $H_i$
in $H = H_1 \amalg \cdots \amalg H_n$, so that we have a correspondence
between $x_i$ and $H_i$. This induces a correspondence between the cluster variables $x'_i$ 
in the cluster $\underline{x}'$ and the indecomposable direct summands 
$T'_i$ in $T'$, which we also denote by $\alpha$. We do not know in general if
the definition of $\alpha$ only depends on the seed $(\underline{x}',B')$.

We can now formulate the connection between cluster algebras and tilting in cluster categories
implied by our main result.

\begin{thm}
Let the notation be as above, with $(\underline{x},B)$ an initial seed for an acyclic cluster algebra,
and $(T', Q_{T'})$ a tilting seed corresponding to a seed $(\underline{x}', B')$, via
the correspondence $\alpha$, inducing a correspondence $x'_i \leftrightarrow T'_i$
for $x'_i \in \underline{x}'$ and $T'_i$ an indecomposable direct summand of $T'$.
\begin{itemize}
\item[(a)]{For any $i \in \{1,\dots,n \}$ we have a commutative diagram
$$
\xy
\xymatrix{
((\underline{x}',B')(i_1,\dots, i_t)) \ar[d]_{\mu_i} \ar[r]^>>>>>>>{\alpha} & (T', Q_{T'}) \ar[d]_{\delta_i} \\
((\underline{x}'',B'')(i_1,\dots, i_t,i)) \ar[r]^>>>>>{\alpha} & (T'', Q_{T''})
}
\endxy
$$
where $\underline{x}''$ is the cluster obtained from $\underline{x}'$ by replacing
$x'_i \in \underline{x}'$ by $x''_i$, and $T''$ is the tilting object in $\C$ obtained
by exchanging the indecomposable summand $T'_i$ by $T''_i$ where
$T'= \bar{T} \amalg T'_i$ and $T'' = \bar{T} \amalg T''_i$ are non-isomorphic tilting objects.
}
\item[(b)]{Identifying $x'_i$ with $T'_i$ and $x''_i$ with $T''_i$, the multiplication rule
for $x'_i x''_i$ is given by
$$T'_i T''_i = \prod (T'_j)^{a_j} + \prod (T'_k)^{c_k}$$
where $a_j$ and $c_k$ are determined by the minimal respectively right and left 
$\add \bar{T}$-approximations $\amalg (T'_j)^{a_j} \to T'_i$ and
$T'_i \to \amalg (T'_k)^{c_k}$.
}
\end{itemize}
\end{thm}

\begin{proof}
(a): This follows by induction, using Theorem \ref{mutate}, where
$\delta_i$ is interpreted as given by a mutation rule like $\mu_i$.

\noindent(b) Let $T'_i$ be the direct summand of $T'$ corresponding
to $x'_i$. By (a), $Q_{T'}$ is the quiver of $B'$, and the monomials 
$M_1$ and $M_2$ are given by the entries of the matrix $B'$,
hence by the arrows in the quiver $Q_{T'}$.
In particular, the arrows entering and leaving $i$, are
given by the minimal right and minimal left 
$\add \bar{T}$-approximations of $T'_i$.
\end{proof}

Note that with the appropriate formulation, this solves Conjecture 9.3 in
\cite{bmrrt}.

For algebras of finite type we know from \cite{bmrrt} that the map $\alpha$
gives a one-one correspondence between the seeds and tilting seeds,
in particular it does not depend on the the $t$-tuple $(i_1, \dots,i_t)$.
In fact, we have in this case a 1--1 correspondence between cluster
variables and indecomposable objects of $\C$, inducing a 
1--1 correspondence between clusters and tilting objects.

Two cluster variables $x_i$ and $x_i^{\ast}$ are said to form an exchange pair 
if there are $n-1$ cluster variables $\{y_1, \dots, y_{n-1} \}$ such that
$\{x_i, y_1, \dots, y_{n-1} \}$ and $\{x_i^{\ast}, y_1, \dots, y_{n-1} \}$ are clusters.
Similarly we have exchange pairs with respect to tilting objects. If $\alpha$
identifies $x_i$ and $x_i^{\ast}$ with $T_i$ and $T_i^{\ast}$, respectively, we then have the
following.

\begin{thm}
For a cluster algebra of finite type, let $\alpha$ be the above correspondence between 
seeds and tilting seeds, and between cluster variables and indecomposable objects in the cluster
category.
\begin{itemize}
\item[(a)]{For any $i \in \{1,\dots,n \}$ we have a commutative diagram
$$
\xy
\xymatrix{
(\underline{x}',B') \ar[d]_{\mu_i} \ar[r]^{\alpha} & (T', Q_{T'}) \ar[d]_{\delta_i} \\
(\underline{x}'',B'') \ar[r]^{\alpha} & (T'', Q_{T''})
}
\endxy
$$
}
\item[(b)]{Identify the cluster variables with the indecomposable objects in $\C$
via $\alpha$. We have $$T_i T_i^{\ast} = \prod (T_j)^{a_j} + \prod (T^{\ast}_k)^{c_k}$$
for an exchange pair $T_i$ and $T_i^{\ast}$ where the $a_j$ and $c_k$
appear in the unique non-split triangles
$$T_i^{\ast} \to \amalg T_j^{a_j} \to T_i \to ,$$
and
$$T_i \to \amalg T_k^{c_k} \to T_i^{\ast} \to $$
in $\C$.
}
\end{itemize}
\end{thm}

\end{document}